\documentclass{article}
\usepackage{amsmath, amssymb}
\usepackage[all]{xy}
\newtheorem{thm}{Theorem}[section]
\newtheorem{cor}[thm]{Corollary}
\newtheorem{prop}[thm]{Proposition}
\newenvironment{dfn}{\medskip\refstepcounter{thm}
\noindent{\bf Definition \thesection.\arabic{thm}\ }}{\medskip}
\newenvironment{proof}[1][,]{\medskip\ifcat,#1
\noindent{{\it Proof}:\ }\else\noindent{\it Proof of #1.\ }\fi}
{\hfill$\square$\medskip}
\newenvironment{remark}[1][Remark]{\begin{trivlist}
\item[\hskip \labelsep {\bfseries #1}]}{\end{trivlist}}

\newenvironment{note}[1][Note]{\begin{trivlist}
\item[\hskip \labelsep {\bfseries #1}]}{\end{trivlist}}
\newenvironment{notes}[1][Notes]{\begin{trivlist}
\item[\hskip \labelsep {\bfseries #1}]}{\end{trivlist}}
\newenvironment{ack}[1][Acknowledgements]{\begin{trivlist}
\item[\hskip \labelsep {\bfseries #1}]}{\end{trivlist}}
\newenvironment{rlist}{\begin{list}{$({\rm \roman{enumi}})$}
{\usecounter{enumi} \setlength{\rightmargin}{10pt}
\setlength{\leftmargin}{40pt} \setlength{\itemsep}{2pt}
\setlength{\parsep}{0pt} \setlength{\labelwidth}{40pt}}}{\end{list}}

\def\eq#1{{\rm(\ref{#1})}}

\def\N{{\mathbb N}}

\def\R{{\mathbb R}}

\def\GL{\mathbin{\rm GL}}

\def\G2{\mathop\textrm{G}_2}

\def\d{{\rm d}}
\def\w{\wedge}
\def\C{{\mathbb C}}
\def\Re{\mathop{\rm Re}\nolimits}

\def\Nhat{\hat{N}}
\def\Uhat{\hat{U}}
\def\That{\hat{T}}

\begin{document}

\title{Coassociative 4-folds with Conical Singularities}
\author{\textsc{Jason Dean Lotay}\\ University College\\ Oxford}
\date{}

\maketitle

\section{Introduction}

This paper is dedicated to the study of \emph{deformations} of
\emph{coassociative 4-folds} in a $\G2$ manifold which have
\emph{conical singularities}.  Understanding the deformations of
such singular coassociative 4-folds should be a useful step towards
attempting to prove a 7-dimensional analogue of the SYZ conjecture.
The research detailed here is motivated by the work on the
deformation theory of special Lagrangian $m$-folds with conical
singularities by Joyce in the series of papers \cite{Joyce1},
\cite{Joyce2}, \cite{Joyce3}, \cite{Joyce4} and \cite{Joyce5}, and
the work of the author in \cite{Lotay} on deformations of
asymptotically conical coassociative 4-folds.

We begin, in Section \ref{sings1}, by discussing the notions of
$\G2$ structures, $\G2$ manifolds and coassociative 4-folds. In
Section \ref{defns} we introduce a distinguished class of singular
manifolds known as \emph{CS manifolds}. CS manifolds have conical
singularities and their nonsingular part is a \emph{noncompact}
Riemannian manifold.  We also define what we mean by CS
coassociative 4-folds.

In order that we may employ various analytic techniques in the
course of our study, we choose to use \emph{weighted Banach spaces}
of forms on the nonsingular part of a CS manifold. These spaces are
described in $\S$\ref{Banach}.  We then focus, in Section
\ref{exceptional}, on a particular linear, elliptic, first-order
differential operator acting between weighted Banach spaces in the
case of a 4-dimensional CS manifold.  The Fredholm and index theory
of this operator is discussed using the theory developed in
\cite{LockhartMcOwen}.

In Section \ref{ch8s2} we stratify the types of deformations allowed
into three problems, each with an associated nonlinear first-order
differential operator whose kernel gives a local description of the
moduli space. The main result for each problem, given in
$\S$\ref{ch8s3}, states that the moduli space is locally
homeomorphic to the kernel of a smooth map between smooth manifolds.
In each case, the map in question can be considered as a projection
from the \emph{infinitesimal deformation space} onto the
\emph{obstruction space}.  Thus, when there are no obstructions the
moduli space is a smooth manifold.  Furthermore, using the material
in $\S$\ref{exceptional} helps to provide a lower bound on the
expected dimension of the moduli space.

The last section shows that, in weakening the condition on the $\G2$
structure of the ambient 7-manifold, there is a generic smoothness
result for the moduli spaces of deformations corresponding to our
second and third problems.

\begin{notes}\begin{itemize}\item[]
\item[(a)] Manifolds are taken to be nonsingular and submanifolds to
be embedded, for convenience, unless stated otherwise.
\item[(b)] We use the convention that the natural numbers
$\N=\{0,1,2,\ldots\}$.\end{itemize} \end{notes}

\section{Coassociative 4-folds}\label{sings1}

The key to defining coassociative 4-folds lies with the introduction
of a distinguished 3-form on $\R^7$.

\begin{dfn}\label{ch2s3dfn1} Let $(x_1,\ldots,x_7)$ be coordinates on $\R^7$ and write
$d{\bf x}_{ij\ldots k}$ for the form $dx_i\w dx_j\w\ldots\w dx_k$.
Define a 3-form $\varphi_0$ by:
\begin{equation}\label{ch2s3eq1}
\varphi_0 = d{\bf x}_{123}+d{\bf x}_{145}+d{\bf x}_{167}+d{\bf
x}_{246}- d{\bf x}_{257}-d{\bf x}_{347}-d{\bf x}_{356}.
\end{equation}
The 4-form $\ast\varphi_0$, where $\varphi_0$ and $\ast\varphi_0$
are related by the Hodge star, is given by:
\begin{equation}\label{ch2s3eq2}
\ast\varphi_0 = d{\bf x}_{4567}+d{\bf x}_{2367}+d{\bf
x}_{2345}+d{\bf x}_{1357}-d{\bf x}_{1346}-d{\bf x}_{1256}-d{\bf
x}_{1247}.
\end{equation}
\end{dfn}

\vspace{-20pt}

\noindent Our choice of expression \eq{ch2s3eq1} for $\varphi_0$
follows that of \cite[Chapter 10]{Joy1}.  This form is sometimes
known as the $\G2$ 3-form because the Lie group $\G2$ is the
subgroup of $\GL(7,\R)$ preserving $\varphi_0$.

\begin{dfn}\label{ch2s3prop3} A 4-dimensional submanifold $N$ of\/ $\R^7$
is coassociative if and only if $\varphi_0|_N\equiv 0$ and
$*\varphi_0|_N>0$.
\end{dfn}

\noindent This definition is not standard but is equivalent to the
usual definition in the language of \emph{calibrated geometry} by
\cite[Proposition IV.4.5 \& Theorem IV.4.6]{HarLaw}.

\begin{remark}
The condition $\varphi_0|_N\equiv 0$ forces $*\varphi_0$ to be a
nonvanishing 4-form on $N$.  Thus, the positivity of $*\varphi_0|_N$
is equivalent to a choice of orientation on $N$.
\end{remark}

So that we may describe coassociative submanifolds of more general
7-manifolds, we make two definitions following \cite[p. 7]{Bryant4}
and \cite[p. 243]{Joy1}.

\begin{dfn}\label{ch2s3subs2dfn1}
Let $M$ be an oriented 7-manifold.  For each $x\in M$ there exists
an orientation preserving isomorphism $\iota_x:T_xM\rightarrow\R^7$.
Since $\text{dim}\,\G2=14$, $\text{dim}\,\GL_+(T_xM)=49$ and
$\text{dim}\,\Lambda^3T^*_xM=35$, the $\GL_+(T_xM)$ orbit of
$\iota_x^*(\varphi_0)$ in $\Lambda^3T^*_xM$, denoted
$\Lambda^3_+T^*_xM$, is open.  A 3-form $\varphi$ on $M$ is
\emph{definite}, or \emph{positive}, if
$\varphi|_{T_xM}\in\Lambda^3_+T^*_xM$ for all $x\in M$.  Denote the
bundle of definite 3-forms $\Lambda^3_+T^*M$.  It is a bundle with
fibre $\GL_+(7,\R)/\G2$ which is \emph{not} a vector subbundle of
$\Lambda^3T^*M$.
\end{dfn}

\noindent Essentially, a definite 3-form is identified with the
$\G2$ 3-form on $\R^7$ at each point in $M$.  Therefore, to each
definite 3-form $\varphi$ we can uniquely associate a 4-form
$*\varphi$ and a metric $g$ on $M$ such that the triple
$(\varphi,*\varphi,g)$ corresponds to $(\varphi_0,*\varphi_0,g_0)$
at each point.  This leads us to our next definition.

\begin{dfn}\label{ch2s3subs2dfn2}
Let $M$ be an oriented 7-manifold, let $\varphi$ be a definite
3-form on $M$ and let $g$ be the metric associated to $\varphi$.  We
call $(\varphi,g)$ a $\G2$ \emph{structure} on $M$. If $\varphi$ is
closed (or coclosed) then $(\varphi,g)$ is a \emph{closed} (or
\emph{coclosed}) $\G2$ structure.  A closed and coclosed $\G2$
structure is called \emph{torsion-free}.
\end{dfn}

\noindent Our choice of notation here agrees with \cite{Bryant4}.

\begin{remark} There is a 1-1
correspondence between pairs $(\varphi,g)$ and principal $\G2$
subbundles of the frame bundle.
\end{remark}

Our definition of torsion-free $\G2$ structure is not standard, but
agrees with other definitions by the following result
\cite[Lemma
11.5]{Salamon}.

\begin{prop}\label{ch2s3subs2prop1} Let $(\varphi,g)$ be a
$\text{\emph{G}}_2$ structure and let $\nabla$ be the Levi--Civita
connection of $g$.  The following are equivalent:
$$d\varphi=d^*\varphi=0;\quad\nabla\varphi=0;\!\quad\!\text{and}\!\quad\!
\text{\emph{Hol}}(g)\subseteq\text{\emph{G}}_2\;\text{with $\varphi$
as the associated 3-form}.$$
\end{prop}

\begin{dfn}\label{ch2s3subs2dfn3} Let $M$ be an oriented
7-manifold endowed with a $\G2$ structure $(\varphi,g)$, denoted
$(M,\varphi,g)$. We say that $(M,\varphi,g)$ is a
\emph{$\varphi$-closed}, or \emph{$\varphi$-coclosed}, 7-manifold if
$(\varphi,g)$ is a closed, respectively coclosed, $\G2$ structure.
If $(\varphi,g)$ is torsion-free, we call $(M,\varphi,g)$ a $\G2$
\emph{manifold}.
\end{dfn}

We are now able to complete our definitions.

\begin{dfn}\label{ch2s3subs2prop2}
A 4-dimensional submanifold $N$ of\/ $(M,\varphi,g)$ is
coassociative if and only if\/ $\varphi|_N\equiv 0$ and
$*\varphi|_N>0$.
\end{dfn}

We end this section with a result, which follows from
\cite[Proposition 4.2]{McLean}, that is invaluable in describing the
deformation theory of coassociative 4-folds.

\begin{prop}\label{ch2s3subs2thm1}
Let $N$ be a coassociative 4-fold in $(M,\varphi,g)$. There is an
isomorphism between the normal bundle $\nu(N)$ of $N$ in $M$ and
$\Lambda^2_+T^*N$ given by $v\mapsto (v\cdot\varphi)|_{TN}$.
\end{prop}

\section{Conical singularities}\label{defns}

\subsection{CS manifolds}

\begin{dfn}\label{ch6s1dfn3}\label{singdfn1} Let $M$ be a connected Hausdorff topological
space and let $z_1,\ldots,z_s\in M$.  Suppose that
$\hat{M}=M\setminus\{z_1,\ldots,z_s\}$ has the structure of a
(nonsingular) $n$-dimensional Riemannian manifold, with Riemannian
metric $g$, compatible with its topology. Then $M$ is a
\emph{
manifold with conical singularities} (at $z_1,\ldots,z_s$ with rate
$\lambda$) if there exist
constants $\epsilon>0$ and $\lambda>1$,
a compact 
$(n\!-\!1)$-dimensional Riemannian manifold $(\Sigma_i,h_i)$,
an open set $U_i\ni z_i$ in $M$ with
$U_i\cap U_j=\emptyset$ for $j\neq i$ and 
a diffeomorphism $\Psi_i:(0,\epsilon)\times\Sigma_i\rightarrow
U_i\setminus\{z_i\}\subseteq\hat{M}$,
for $i=1,\ldots,s$, such that
 \begin{equation}\label{ch6s1eq2}
|\nabla_i^j(\Psi_i^*(g)-g_i)|=O(r_i^{\lambda-1-j})\qquad \text{for
$j\in\N$ as $r_i\rightarrow 0$,}
\end{equation}
 where $r_i$ is the coordinate on $(0,\infty)$ on
the cone $C_i=(0,\infty)\times\Sigma_i$, $g_i=dr_i^2+r_i^2h_i$ is
the conical metric on $C_i$, $\nabla_i$ is the Levi--Civita
connection derived from $g_i$ and $|.|$ is calculated using $g_i$.
 We call $C_i$ the \emph{cone} at the singularity $z_i$ and
let the \emph{ends} $\hat{M}_{\infty}$ of $\hat{M}$ be the disjoint
union
$$\hat{M}_{\infty}=\bigsqcup_{i=1}^s U_i\setminus\{z_i\}.$$

We say that $M$ is \emph{CS} or a \emph{CS manifold} (with rate
$\lambda$) if it is a
manifold with conical singularities which have rate $\lambda$ and it
is compact as a topological space.  In these circumstances it may be
written as the disjoint union
$$M=K\sqcup \bigsqcup_{i=1}^sU_i,$$
where $K$ is compact as it is closed in $M$.
\end{dfn}

The condition $\lambda>1$ guarantees that the metric on $\hat{M}$
genuinely converges to the conical metric on $C_i$, as is evident
from \eq{ch6s1eq2}.  Since $M$ is supposed to be Hausdorff, the set
$U_i\setminus\{z_i\}$ is open in $\hat{M}$ for all $i$. Moreover,
the condition that the $U_i$ are disjoint may be easily satisfied
since, if $i\neq j$, $z_i$ and $z_j$ may be separated by two
disjoint open sets and, by hypothesis, there are only a finite
number of singularities.

\begin{remark} If $M$ is a CS manifold, $\hat{M}$ is a
\emph{noncompact} manifold.
\end{remark}

\begin{dfn}\label{ch6s1dfn4}\label{radiusfn}
Let $M$ be a CS manifold. Using the notation of Definition
\ref{ch6s1dfn3}, a \emph{radius function} on $\hat{M}$ is a smooth
function $\rho:\hat{M}\rightarrow (0,1]$, bounded below by a
positive constant on $\hat{M}\setminus\hat{M}_{\infty}$, such that
there exist positive constants $c_1<1$ and $c_2>1$ with
$$c_1r_i<\Psi_i^*(\rho)<c_2r_i$$
on $(0,\epsilon)\times\Sigma_i$ for $i=1,\ldots,s$.
\end{dfn}

If $M$ is CS we may construct a radius function on $\hat{M}$ as
follows. Let $\rho(x)=1$ for all
$x\in\hat{M}\setminus\hat{M}_{\infty}$. Define
$\rho_i:\Psi_i((0,\epsilon/2)\times\Sigma_i)\rightarrow (0,1)$ to be
equal to $r_i/\epsilon$ for $i=1,\ldots,s$ and then define $\rho$ by
interpolating smoothly between its definition on
$\hat{M}\setminus\hat{M}_{\infty}$ and $\rho_i$ on each of the
disjoint sets $\Psi_i((\epsilon/2,\epsilon)\times\Sigma_i)$.

\subsection{CS coassociative 4-folds}

Let $B(0;\eta)$ denote the open ball about $0$ in $\R^7$ with radius
$\eta>0$, i.e. $B(0;\eta)=\{\textbf{v}\in\R^7:|\textbf{v}|<\eta\}$.
We define a preferred choice of local coordinates on a $\text{G}_2$
manifold near a finite set of points.

\begin{dfn}\label{ch8s1dfn1}\label{coords}
Let $(M,\varphi,g)$ be a $\text{G}_2$ manifold as in Definition
\ref{ch2s3subs2dfn3} and let $z_1,\ldots,z_s$ be points in $M$.
There exist
a constant $\eta>0$, 
an open set
$V_i\ni z_i$ in $M$ with $V_i\cap V_j=\emptyset$ for $j\neq i$ and 
a diffeomorphism $\chi_i:B(0;\eta)\subseteq\R^7\rightarrow V_i$ with
$\chi_i(0)=z_i$,
for $i=1,\ldots,s$, such that $\zeta_i=d\chi_i|_0:\R^7\rightarrow
T_{z_i}M$ is an isomorphism identifying the standard $\text{G}_2$
structure $(\varphi_0,g_0)$ on $\R^7$ with the pair
$(\varphi|_{T_{z_i}M},g|_{T_{z_i}M})$.
We call the set $\{\chi_i:B(0;\eta)\rightarrow V_i:i=1,\ldots,s\}$ a
\emph{$\text{\emph{G}}_2$ coordinate system near
${z_1,\ldots,z_s}$}.

We say that two $\text{G}_2$ coordinate systems near
$z_1,\ldots,z_s$, with maps $\chi_i$ and $\tilde{\chi}_i$ for
$i=1,\ldots,s$ respectively, are \emph{equivalent} if
$d\tilde{\chi}_i|_0=d\chi_i|_0=\zeta_i$ for all $i$.
\end{dfn}

\noindent The definition above is an analogue of the local
coordinate system for almost Calabi--Yau manifolds used by Joyce
\cite[Definition 3.6]{Joyce1}.  Although the family of $\text{G}_2$
coordinate systems near $z_1,\ldots,z_s$ is clearly
infinite-dimensional, there are only finitely many equivalence
classes, given by the number of possible sets
$\{\zeta_1,\ldots,\zeta_s\}$. Moreover, the family of choices for
each $\zeta_i$ is isomorphic to $\text{G}_2$. \begin{note}
Definition \ref{ch8s1dfn1} does not require the $\G2$ structure
$(\varphi,g)$ to be \emph{torsion-free}.\end{note}

\begin{dfn}\label{ch8s1dfn2}\label{singdfn2}
Let $(M,\varphi,g)$ be a $\text{G}_2$ manifold, let $N\subseteq M$
be compact and connected and let $z_1,\ldots,z_s\in N$. We say that
$N$ is a 4-fold in $M$ with \emph{conical singularities at
$z_1,\ldots,z_s$
with rate $\lambda$}, denoted a \emph{CS 4-fold}, if 
$\hat{N}=N\setminus\{z_1,\ldots,z_s\}$ is a (nonsingular)
4-dimensional submanifold of $M$ and there exist
constants $0<\epsilon<\eta$ and
$\lambda>1$, 
a compact 3-dimensional Riemannian submanifold $(\Sigma_i,h_i)$ of
$\mathcal{S}^6\subseteq\R^7$, where $h_i$ is the restriction of the
round metric on $\mathcal{S}^6$ to $\Sigma_i$,
an open set $U_i\ni z_i$ in $N$ with $U_i\subseteq
V_i$ 
and
a smooth map $\Phi_i:(0,\epsilon)\times\Sigma_i\rightarrow
B(0;\eta)\subseteq\R^7$,
for $i=1,\ldots,s$,
 such that
 $\Psi_{i}=\chi_i\circ\Phi_i:(0,\epsilon)\times\Sigma_i\rightarrow
 U_i\setminus\{z_i\}$ is a diffeomorphism and
 $\Phi_i$ satisfies
\begin{equation}\label{ch8s1eq1}
|\nabla^j_i(\Phi_i(r_i,\sigma_i)-\iota_i(r_i,\sigma_i))|
=O(r_i^{\lambda-j})\qquad \text{for $j\in\N$ as $r_i\rightarrow 0$,}
\end{equation}
where $\iota_i(r_i,\sigma_i)=r_i\sigma_i\in B(0;\eta)$, $\nabla_i$
is the Levi--Civita connection of the cone metric
$g_i=dr_i^2+r_i^2h_i$ on $C_i=(0,\infty)\times\Sigma_i$ coupled with
partial differentiation on $\R^7$,
 $|.|$ is calculated with respect to $g_i$
  and $\{\chi_i:B(0;\eta)\rightarrow
V_i:i=1,\ldots,s\}$ is a $\text{G}_2$ coordinate system near
$z_1,\ldots,z_s$.

 We call $C_i$ the \emph{cone} at the singularity $z_i$ and
$\Sigma_i$ the \emph{link} of the cone $C_i$.  We may write $N$ as
the disjoint union
$$N=K\sqcup\bigsqcup_{i=1}^sU_i ,$$
where $K$ is compact.

If $\Nhat$ is coassociative in $M$, we say that $N$ is a \emph{CS
coassociative 4-fold}.
\end{dfn}

Suppose $N$ is a CS 4-fold at $z_1,\ldots,z_s$ with rate $\lambda$
in $(M,\varphi,g)$ and use the notation of Definition
\ref{ch8s1dfn2}. The induced metric on $\Nhat$, $g|_{\Nhat}$, makes
$\Nhat$ into a Riemannian manifold. Moreover, it is clear from
\eq{ch8s1eq1} that the maps $\Psi_i$ satisfy \eq{ch6s1eq2} in
Definition \ref{ch6s1dfn3} with the same constant $\lambda$. Thus,
$N$ may be considered as a CS manifold with rate $\lambda$.

It is important to note that, if $\lambda\in(1,2)$, Definition
\ref{ch8s1dfn2} is independent of the choice of $\text{G}_2$
coordinate system near the singularities, up to equivalence. Suppose
we have two equivalent coordinate systems defined using maps
$\chi_i$ and $\tilde{\chi}_i$.  These maps must agree up to second
order since the zero and first order behaviour of each is
prescribed, as stated in Definition \ref{ch8s1dfn1}. Therefore, the
transformed maps $\tilde{\Phi}_i$ corresponding to $\tilde{\chi}_i$
such that $\tilde{\Psi}_i=\tilde{\chi}_i\circ\tilde{\Phi}_i=
\chi_i\circ\Phi_i=\Psi_i$ are defined by:
$$\tilde{\Phi}_i=(\tilde{\chi}_i^{-1}\circ\chi_i)\circ\Phi_i.$$
Hence
$$|\nabla^j_i(\tilde{\Phi}_i(r_i,\sigma_i)-\Phi_i(r_i,\sigma_i))|
=O(r_i^{2-j})\qquad\text{for $j\in\N$ as $r_i\rightarrow0$,}$$ where
$\nabla_i$ and $|.|$ are calculated as in Definition
\ref{ch8s1dfn2}.  Thus, in order that the terms generated by the
transformation of the $\text{G}_2$ coordinate system neither
dominate nor be of equal magnitude to the
$O(r_i^{\lambda-j}\hspace{1pt})$ terms given in \eq{ch8s1eq1}, we
need $\lambda<2$.

\medskip

We now make a definition which also depends only on equivalence
classes of $\text{G}_2$ coordinate systems near the singularities.

\begin{dfn}\label{ch8s1dfn3}\label{tangent}
Let $N$ be a CS 4-fold at $z_1,\ldots,z_s$ in a $\text{G}_2$
manifold $(M,\varphi,g)$.  Use the notation of Definitions
\ref{ch8s1dfn1} and \ref{ch8s1dfn2}. For $i=1,\ldots,s$ define a
cone $\hat{C}_i$ in $T_{z_i}M$ by
$\hat{C}_i=(\zeta_i\circ\iota_i)(C_i)$. We call $\hat{C}_i$ the
\emph{tangent cone} at $z_i$.
\end{dfn}

\noindent One can show that $\hat{C}_i$ is a tangent cone to $N$ at
$z_i$ in the sense of \emph{geometric measure theory} (see, for
example, \cite[p. 233]{Federer}).  We also have a straightforward
result relating to the tangent cones at singular points of CS
coassociative 4-folds.

\begin{prop}\label{ch8s1prop1}
Let $N$ be a CS coassociative 4-fold at $z_1,\ldots,z_s$ in a
$\text{\emph{G}}_2$ manifold $(M,\varphi,g)$.  The tangent cones at
$z_1,\ldots,z_s$ are coassociative.
\end{prop}

\begin{proof} Use the notation of Definitions \ref{ch8s1dfn1} and \ref{ch8s1dfn2}.

It is enough to show that $\iota_i(C_i)$ is coassociative in $\R^7$
for all $i$, since $\zeta_i:\R^7\rightarrow T_{z_i}M$ is an
isomorphism identifying $(\varphi_0,g_0)$ with
$(\varphi|_{T_{z_i}M},g|_{T_{z_i}M})$.  This is equivalent to the
condition $\iota_i^*(\varphi_0)\equiv0$ for $i=1,\ldots,s$.

Note that $\varphi|_{\Nhat}\equiv 0$ implies that, for all $i$,
$\varphi|_{U_i\setminus\{z_i\}}\equiv 0$. Hence,
$\Psi_i^*(\varphi)=\Phi_i^*(\chi_i^*(\varphi))$ vanishes on $C_i$
for all $i$. Using \eq{ch8s1eq1},
$$|\Phi_i^*(\chi_i^*(\varphi))-\iota_i^*(\chi_i^*(\varphi))|
=O(r_i^{\lambda-1})\qquad\text{as $r_i\rightarrow0$}$$ for all $i$.
Moreover,
$$|\iota_i^*(\chi_i^*(\varphi))-
\iota_i^*(\varphi_0)|=O(r_i) \qquad\text{as $r_i\rightarrow0$}$$
since $$\chi_i^*(\varphi)=\varphi_0+O(r_i)\quad\text{and}\quad
|\nabla\iota_i|=O(1) \qquad\text{as $r_i\rightarrow0$}.$$ Therefore,
because $\lambda>1$,
$$|\iota_i^*(\varphi_0)|\rightarrow0\qquad\text{as
$r_i\rightarrow 0$}$$ for all $i$.  As
$T_{r_i\sigma_i}\iota_i(C_i)=T_{\sigma_i}\iota_i(C_i)$ for all
$(r_i,\sigma_i)\in C_i$, $|\iota_i^*(\varphi_0)|$ is independent of
$r_i$ and thus vanishes for all $i$ as required.
\end{proof}

\section{Weighted Banach spaces}\label{Banach}

For this section let $M$ be an $n$-dimensional CS manifold and let
$\hat{M}$ be its nonsingular part as in Definition \ref{singdfn1}.
We define \emph{weighted} Banach spaces of forms as in
\cite[$\S$1]{Bartnik}, as well as the usual `unweighted' spaces.

\begin{dfn}\label{ch6s2dfn1}\label{sobolev}
Let $p\geq 1$ and let $k,m\in\N$ with $m\leq n$. The \emph{Sobolev
space} $L_k^p(\Lambda^mT^*\hat{M})$ is the set of $m$-forms $\xi$ on
$\hat{M}$ which are $k$ times weakly differentiable and such that
the norm
\begin{equation}\label{ch6s2eq1}
\|\xi\|_{L_k^p}=\left(\sum_{j=0}^k\int_{\hat{M}} |\nabla^j\xi|^p
\,dV_g\right)^\frac{1}{p} \end{equation} is finite.  The normed
vector space $L_k^p(\Lambda^mT^*\hat{M})$ is a Banach space for all
$p\geq 1$ and $L_k^2(\Lambda^mT^*\hat{M})$ is a Hilbert space.

We introduce the space of $m$-forms
\begin{equation*}
L_{k,\,\text{loc}}^p(\Lambda^mT^*\hat{M})=\{\xi\,:\, f\xi\in
L_k^p(\Lambda^mT^*\hat{M})\;\text{for all} \,f\in
C_{\text{cs}}^{\infty}(\hat{M})\}
\end{equation*}
where $C_{\text{cs}}^{\infty}(\hat{M})$ is the space of smooth
functions on $\hat{M}$ with compact support.

Let $\mu\in\R$ and let $\rho$ be a radius function on $\hat{M}$. The
\emph{weighted Sobolev space} $L_{k,\,\mu}^p(\Lambda^mT^*\hat{M})$
of $m$-forms $\xi$ on $\hat{M}$ is the subspace of
$L^p_{k,\,\text{loc}}(\Lambda^mT^*\hat{M})$ such that the norm
\begin{equation}\label{ch6s2eq2}
\|\xi\|_{L_{k,\,\mu}^p}=\left(\sum_{j=0}^k\int_{\hat{M}}
|\rho^{j-\mu}\nabla^j\xi|^p\rho^{-n} \,dV_g\right)^\frac{1}{p}
\end{equation}
is finite. Then $L_{k,\,\mu}^p(\Lambda^mT^*\hat{M})$ is a Banach
space and $L_{k,\,\mu}^2(\Lambda^mT^*\hat{M})$ is a Hilbert space.
\end{dfn}

We may note here, trivially, that $L_0^p(\Lambda^mT^*\hat{M})$ is
equal to the standard $L^p$-space of $m$-forms on $\hat{M}$.
Further, by comparing equations \eq{ch6s2eq1} and \eq{ch6s2eq2} for
the respective norms,
$L^p(\Lambda^mT^*\hat{M})=L_{0,\,-\frac{n}{p}}^p(\Lambda^mT^*\hat{M})$.
In particular,
\begin{equation}\label{ch6s2eq3}
L^2(\Lambda^mT^*\hat{M})=L_{0,\,-\frac{n}{2}}^2(\Lambda^mT^*\hat{M}).
\end{equation}

For the following two definitions we take
$C^k_{\text{loc}}(\Lambda^mT^*\hat{M})$ to be the vector space of
$k$ times continuously differentiable $m$-forms.

\begin{dfn}\label{ch6s2dfn2} Let $\rho$ be a radius function on
$\hat{M}$, let $\mu\in\R$ and let $k,m\in\N$ with $m\leq n$.  The
\emph{weighted $C^k$-space} $C_{\mu}^{k}(\Lambda^mT^*\hat{M})$ of
$m$-forms $\xi$ on $\hat{M}$ is the subspace of
$C^k_{\text{loc}}(\Lambda^mT^*\hat{M})$ such that the norm
$$\|\xi\|_{C_{\mu}^{k}}=\sum_{j=0}^k
\sup_{\hat{M}}|\rho^{j-\mu}\nabla^j\xi|$$is finite. We also define
$$C_{\mu}^{\infty}(\Lambda^mT^*\hat{M})=\bigcap_{k\geq
0}C_{\mu}^{k}(\Lambda^mT^*\hat{M}).$$  Then
$C_{\mu}^{k}(\Lambda^mT^*\hat{M})$ is a Banach space but in general
$C_{\mu}^{\infty}(\Lambda^mT^*\hat{M})$ is not.
\end{dfn}

In the next definition we refer to the usual normed vector space
$C^k(\Lambda^mT^*\hat{M})$ of $k$ times continuously differentiable
$m$-forms such that the following norm is finite:
$$\|\xi\|_{C^k}= \sum_{j=0}^k \sup_{\hat{M}}|\nabla^j\xi|.$$

\begin{dfn}\label{ch6s2dfn3}
Let $d(x,y)$ be the geodesic distance between points $x,y\in
\hat{M}$ and let $\rho$ be a radius function on $\hat{M}$. Let $a\in
(0,1)$ and let $k,m\in\N$ with $m\leq n$. Let
\begin{align*}H=\{&(x,y)\in \hat{M}\times \hat{M}\,:\,x\neq
y,\,c_1\rho(x)\leq\rho(y)\leq c_2\rho(x)\,\;\text{and}\;\,\\
&\text{there exists a geodesic in $\hat{M}$ of length $d(x,y)$ from
$x$ to $y$}\},\end{align*} where $0<c_1<1<c_2$ are constant. A
section $s$ of a vector bundle $V$ on $\hat{M}$, endowed with a
connection, is \emph{H\"older continuous} (with \emph{exponent $a$})
if
$$[s]^a=
\sup_{(x,y)\in H}\frac{|s(x)-s(y)|_{V}}{d(x,y)^a}<\infty.$$ We
understand the quantity $|s(x)-s(y)|_V$ as follows.  Given $(x,y)\in
H$, there exists a geodesic $\gamma$ of length $d(x,y)$ connecting
$x$ and $y$. Parallel translation along $\gamma$ using the
connection on $V$ identifies the fibres over $x$ and $y$ and the
metrics on them. Thus, with this identification, $|s(x)-s(y)|_V$ is
well-defined.

The \emph{H\"older space} $C^{k,\,a}(\Lambda^mT^*\hat{M})$ is the
set of $\xi\in C^k(\Lambda^mT^*\hat{M})$ such that $\nabla^k\xi$ is
H\"older continuous (with exponent $a$) and the norm
$$\|\xi\|_{C^{k,\,a}}=\|\xi\|_{C^k}+[\nabla^k\xi]^a$$
is finite. The normed vector space $C^{k,\,a}(\Lambda^mT^*\hat{M})$
is a Banach space.

We also introduce the notation
\begin{align*}
C^{k,\,a}_{\text{loc}}&(\Lambda^mT^*\hat{M})\\&=\!\{\xi\in
C^k_{\text{loc}}(\Lambda^mT^*\hat{M}): f\xi\in
C^{k,\,a}(\Lambda^mT^*\hat{M})\;\text{for all} \,f\in
C_{\text{cs}}^{\infty}(\hat{M})\}.\end{align*}

Let $\mu\in\R$.  The \emph{weighted H\"older space}
$C_{\mu}^{k,\,a}(\Lambda^mT^*\hat{M})$ of $m$-forms $\xi$ on
$\hat{M}$ is the subspace of
$C^{k,\,a}_{\text{loc}}(\Lambda^mT^*\hat{M})$ such that the norm
$$\|\xi\|_{C^{k,\,a}_{\mu}}=\|\xi\|_{C^{k}_{\mu}}+[\xi]^{k,\,a}_{\mu}$$
is finite, where
$$[\xi]^{k,\,a}_{\mu}=[\rho^{k+a-\mu}\nabla^k\xi]^{a}.
$$  Then $C_{\mu}^{k,\,a}(\Lambda^mT^*\hat{M})$ is a
Banach space.  It is clear that we have an embedding
$C_{\mu}^{k,\,a}(\Lambda^mT^*\hat{M})\hookrightarrow
C_{\mu}^l(\Lambda^mT^*\hat{M})$ whenever $l\leq k$.
\end{dfn}

We shall need the analogue of the Sobolev Embedding Theorem for
weighted spaces, which is adapted from \cite[Lemma
7.2]{LockhartMcOwen} and \cite[Theorem 1.2]{Bartnik}.
\begin{thm}[\textbf{Weighted Sobolev Embedding Theorem}]
\label{ch6s2thm1}\label{wembed} Let $p,\,q\geq 1$,\\ $a\in (0,1)$,
$\mu,\nu\in\R$ and $k,l,m\in\N$ with $m\leq n$.
\begin{itemize}
\item[{\rm (a)}] If $k\geq l$, $k-\frac{n}{p}\geq l-\frac{n}{q}$
and either
\begin{rlist}
\item $p\leq q$ and $\mu\geq\nu$
 or \item $p>q$ and $\mu>\nu$,
\end{rlist}
there is a continuous embedding
$L_{k,\,\mu}^p(\Lambda^mT^*\hat{M})\hookrightarrow
L_{l,\,\nu}^q(\Lambda^mT^*\hat{M})$. \item[{\rm (b)}] If
$k-\frac{n}{p}\geq l+a$, there is a continuous embedding
$L_{k,\,\mu}^p(\Lambda^mT^*\hat{M})\hookrightarrow
C_{\mu}^{l,\,a}(\Lambda^mT^*\hat{M})$.
\end{itemize}
\end{thm}

We shall also require an Implicit Function Theorem for Banach
spaces, which follows immediately from \cite[Chapter 6, Theorem
2.1]{Lang2}.
\begin{thm}[\textbf{Implicit Function Theorem}]
\label{ch6s2thm2}\label{implicit} $\!\!$Let $X$ and $Y$ be Banach
spaces and let $W\subseteq X$ be an open neighbourhood of $0$. Let
$\mathcal{G}:W\rightarrow Y$ be a $C^k$ map $(k\geq 1)$ such that
$\mathcal{G}(0)=0$. Suppose further that
$d\mathcal{G}|_{0}:X\rightarrow Y$ is surjective with kernel $K$
such that $X=K\oplus A$ for some closed subspace $A$ of $X$. There
exist open sets $V\subseteq K$ and $V^\prime\subseteq A$, both
containing $0$, with $V\times V^\prime\subseteq W$, and a unique
$C^k$ map $\mathcal{V}:V\rightarrow V^\prime$ such that
$$\text{\emph{Ker}}\,\mathcal{G}\cap(V\times V^\prime)=\{(x,\mathcal{V}(x))\,:\,x\in V\}
$$ in $X=K\oplus A$.
\end{thm}

\section{The operator $d+d^*$}\label{exceptional}

In this section we let $M$ be a \emph{4-dimensional} CS manifold and
let $\hat{M}$ be as in Definition \ref{ch6s1dfn3}.  An essential
part of our study is the use of the Fredholm and index theory for
the elliptic operator $d+d^*$ acting from
$\Lambda^2_+T^*\hat{M}\oplus\Lambda^4T^*\hat{M}$ to
$\Lambda^3T^*\hat{M}$. We therefore consider
\begin{equation}\label{ch6s3eq1}
d+d^*:L_{k+1,\,\mu}^{p}(\Lambda_+^2T^*\hat{M}\oplus
\Lambda^4T^*\hat{M})\rightarrow
L_{k,\,\mu-1}^{p}(\Lambda^3T^*\hat{M}),
\end{equation}
where $p\geq 2$, $k\in\N$ and $\mu\in\R$.

\subsection{Fredholm theory}\label{exceptsubs1}

Our first result follows from \cite[Theorem 1.1 \& Theorem
6.1]{LockhartMcOwen}.
\begin{prop}\label{ch6s3thm1}
There exists a countable discrete set $\mathcal{D}\subseteq\R$ such
that \eq{ch6s3eq1} is Fredholm if and only if
$\mu\notin\mathcal{D}$.
\end{prop}
Moreover, we can give an explicit description of $\mathcal{D}$ by a
similar argument to \cite[p. 13-14]{Lotay}, which is for
\emph{asymptotically conical} (AC) manifolds, as follows.

Recall the notation of Definition \ref{ch6s1dfn3}.  Transform the
metric on $\hat{M}$ to a conformally equivalent metric which is
asymptotically \textit{cylindrical} on the ends $\hat{M}_\infty$ of
$\hat{M}$; that is, if $(t_i,\sigma_i)$ are coordinates on
$(0,\infty)\times\Sigma_i$, the metric is asymptotic to
$dt_i^2+h_i$. With respect to this new metric, $d+d^*$ corresponds
to
$$(d+d^*)_{\infty}=e^{mt}(d+e^{-2t}d^*)e^{-mt}$$
acting on $m$-forms on $\hat{M}$.

Let $$\Sigma=\bigsqcup_{i=1}^s\Sigma_i.$$ If
$\pi:(0,\infty)\times\Sigma\rightarrow\Sigma$ is the natural
projection map, the action of $(d+d^*)_{\infty}$ on
$\pi^*(\Lambda^2T^*\Sigma)\oplus\pi^*(\Lambda^{\text{odd}}T^*\Sigma)$
is:
\begin{equation}\label{ch6s3eq2}
(d+d^*)_{\infty}=\left(\begin{array}{cc}\displaystyle d+d^*&\frac{\partial}{\partial t}+3-m \\
-(\frac{\partial}{\partial t}+m) & -(d+d^*)\end{array}\right)
\end{equation}
where $m$ denotes the operator which multiplies $m$-forms by a
factor $m$.  However, we wish only to consider elements of
$\Lambda^1T^*\Sigma\oplus\Lambda^2T^*\Sigma$ which correspond to
self-dual 2-forms on $\hat{M}$, so we define
$V_{\Sigma}\subseteq\Lambda^2T^*\Sigma\oplus\Lambda^{\text{odd}}T^*\Sigma$
by
\begin{equation*}\label{ch6s3eq3}
V_{\Sigma}=\bigsqcup_{i=1}^s\{(\alpha,\ast\alpha+\beta)\,:\,\alpha\in\Lambda^2T^*\Sigma_i,\,
\beta\in\Lambda^3T^*\Sigma_i\}.
\end{equation*}
Then $\pi^*(V_{\Sigma})$ corresponds to
$\Lambda_+^2T^*\hat{M}\oplus\Lambda^4T^*\hat{M}$.

For $w\in\C$ define a map $(d+d^*)_{\infty}(w)$ by:
\begin{equation}\label{ch6s3eq4}
(d+d^*)_{\infty}(w)=\left(\begin{array}{cc}\displaystyle d+d^*&-w+3-m \\
w-m & -(d+d^*)\end{array}\right)
\end{equation}
 acting on
$V_{\Sigma}\otimes\C$.
 Notice that we have
formally substituted $w$ for $-\frac{\partial}{\partial t}$ in
\eq{ch6s3eq2}.

Let
\begin{equation*}
W_{\Sigma}=\bigsqcup_{i=1}^s\{(*\alpha+\beta,\alpha)\,:\,\alpha\in
\Lambda^2T^*\Sigma_i,\,\beta\in\Lambda^3T^*\Sigma_i\}\subseteq
\Lambda^{\text{odd}}T^*\Sigma\oplus\Lambda^{2}T^*\Sigma.
\end{equation*}
Define $\mathcal{C}\subseteq\C$ as the set of $w$ for which the map
$$
(d+d^*)_{\infty}(w):L_{k+1}^p(V_{\Sigma}\otimes\C) \rightarrow L_k^p
(W_{\Sigma}\otimes\C)
$$ is \emph{not} an isomorphism. By the proof of \cite[Theorem
1.1]{LockhartMcOwen}, $\mathcal{D}=\{\Re w:w\in\mathcal{C}\}$. By
\cite[Lemma 6.1.13]{Marshall}, the corresponding sets
$\mathcal{C}(\Delta^m)$, where $\Delta^m$ is the Laplacian on
$m$-forms, are all real for an asymptotically conical manifold.
Since the same will be true for the CS case, we deduce that
$\mathcal{C}\subseteq\R$. Hence $\mathcal{C}=\mathcal{D}$.

The symbol, hence the index $\text{ind}_w$, of $(d+d^*)_{\infty}(w)$
is independent of $w$. Furthermore, $(d+d^*)_{\infty}(w)$ is an
isomorphism for generic values of $w$ since $\mathcal{D}$ is
countable and discrete. Therefore $\text{ind}_w=0$ for all $w\in\C$;
that is,
\begin{equation*}
\text{dim Ker}(d+d^*)_\infty(w)=\text{dim Coker}(d+d^*)_\infty(w),
\end{equation*}
so that \eq{ch6s3eq4} is not an isomorphism precisely when it is not
injective.

The condition $(d+d^*)_{\infty}(w)=0$, using \eq{ch6s3eq4},
corresponds to the existence of $\alpha\in
C^{\infty}(\Lambda^2T^*\Sigma_i)$ and $\beta\in
C^{\infty}(\Lambda^3T^*\Sigma_i)$, for some $i$, satisfying
\begin{align}
\label{ch6s3eq5}
d\alpha=w\beta\quad\text{and}\quad
d\!*\!\alpha+d^*\beta=(w-2)\alpha.
\end{align}

\begin{notes}\begin{itemize}\item[] \item[(a)] The equations above imply that
$$dd^*\beta=\Delta\beta=w(w-2)\beta.$$ Since eigenvalues of the
Laplacian on $\Sigma_i$ must necessarily be positive, $\beta=0$ if
$w\in(0,2)$.
\item[(b)] If $w=0$ and we take $\alpha=0$,
\eq{ch6s3eq5} forces $\beta$ to be coclosed. As there are nontrivial
coclosed 3-forms on $\Sigma_i$, $(d+d^*)_{\infty}(0)$ is not
injective, so $0\in\mathcal{D}$.
\item[(c)] Suppose that
$w=2$ lies in $\mathcal{D}$. Then \eq{ch6s3eq5} gives $[\beta]=0$ in
$H_{\text{dR}}^3(\Sigma_i)$. We know that $\beta$ is harmonic so, by
Hodge theory, $\beta=0$. Therefore $2\in\mathcal{D}$ if and only if
there exists a nonzero closed and coclosed 2-form on $\Sigma_i$ for
some $i$.\end{itemize}
\end{notes}

We state a proposition which follows from the work above.

\begin{prop}\label{ch6s3prop1}
Let $M$ be a 4-dimensional CS manifold.  Use the notation of
Definition \ref{ch6s1dfn3}.  For $i=1,\ldots,s$ let
$D(\mu,i)=\{(\alpha,\beta)\in C^\infty(\Lambda^2T^*\Sigma_i\oplus
\Lambda^3T^*\Sigma_i)\,:\,d\alpha=\mu\beta,\;d\!*\!\alpha+d^*\beta=
(\mu-2)\alpha\}.$ The set $\mathcal{D}$ of real numbers $\mu$ such
that \eq{ch6s3eq1} is not Fredholm is given by:
$$\mathcal{D}=\bigcup_{i=1}^s\{\mu\in\R\,:\,D(\mu,i)\neq 0\}.$$
\end{prop}

\begin{remark} A perhaps more illuminating way to characterise $D(\mu,i)$ is
by:
\begin{align*}
 (\alpha,\beta)\in D(\mu,i) \Longleftrightarrow&
\;\xi=(r^{\mu-2}\alpha+r^{\mu-1}dr\wedge*\alpha,
r^{\mu-3}dr\wedge\beta)\\&\;\text{is an $O(r^{\mu})$ solution of
$(d+d^*)\xi=0$ on $C_i$,}
\end{align*}
using the notation of Definition \ref{ch6s1dfn3}.
\end{remark}

Lockhart and McOwen \cite[$\S$10]{LockhartMcOwen} study the
Laplacian on $m$-forms on a manifold with a conical singularity.
From this work, which can easily be extended to manifolds with more
than one singularity, we can make an important observation about the
set $\mathcal{D}$.

\begin{prop}\label{nopoints}In the notation of Proposition
\ref{ch6s3prop1}, $\mathcal{D}\cap(-2,-1]=\emptyset$.
\end{prop}

\begin{proof}
Let
$$\Delta^m:L^p_{k+1,\,\mu}(\Lambda^mT^*\hat{M})\longrightarrow
L^p_{k-1,\,\mu-2}(\Lambda^mT^*\hat{M})$$ be the Laplacian on
$m$-forms and denote the set of $\mu$ such that it is not Fredholm
by $\mathcal{D}(\Delta^m)$.  Since $\mu>-2$ and $p\geq 2$ we see
that $L^p_{k+1,\,\mu}\hookrightarrow L^2_{0,\,-2}=L^2$ by Theorem
\ref{ch6s2thm1} and \eq{ch6s2eq3}.

We then apply \cite[Theorem 10.2]{LockhartMcOwen} for the Laplacian
on 2-forms and 4-forms on a 4-dimensional CS manifold to see that
$$\mathcal{D}(\Delta^2)\cap(-2,-1]=\mathcal{D}(\Delta^4)\cap(-2,-1]=\emptyset.$$
Note that our rate $\mu$ is related to the weighting factor in
\cite[$\S$10]{LockhartMcOwen}, which we may denote as $\nu$, by
$\mu=-\nu-2$.  As it is clear that
$\mathcal{D}\subseteq(\mathcal{D}(\Delta^2)\cup\mathcal{D}(\Delta^4))$,
the result follows.
\end{proof}

\subsection{Index theory}\label{indextheory}

We begin with some definitions following \cite{LockhartMcOwen}.

\begin{dfn}\label{ch6s3dfn1}
Use the notation of $\S$\ref{exceptsubs1}. Let $\mu\in\mathcal{D}$.
Define $\d(\mu)$ to be the dimension of the vector space of
solutions of $(d+d^*)_{\infty}\xi=0$ of the form
$$\xi(t,\sigma)=e^{-\mu t}p\,(t,\sigma)$$
where $p\,(t,\sigma)$ is a polynomial in $t\in(0,\infty)$ with
coefficients in $C^{\infty}(V_{\Sigma}\otimes\C)$. 
\end{dfn}

The next result is immediate from \cite[Theorem
1.2]{LockhartMcOwen}.
\begin{thm}\label{ch6s3thm2}
Let $\lambda,\lambda^{\prime}\notin\mathcal{D}$ with
$\lambda^{\prime}\leq\lambda$, where $\mathcal{D}$ is given in
Proposition \ref{ch6s3prop1}. For any $\mu\notin\mathcal{D}$ let
$\text{\emph{ind}}_{\mu}(d+d^*)$ denote the Fredholm index of
\eq{ch6s3eq1}.  Then
$$\text{\emph{ind}}_{\lambda^{\prime}}(d+d^*)-\text{\emph{ind}}_{\lambda}(d+d^*)=\sum_{\mu\in\mathcal{D}\,
\cap(\lambda^{\prime},\,\lambda)}\!\!\!\!\!\!\d(\mu).$$
\end{thm}

We make a key observation, which shall be used on a number of
occasions in later sections.

\begin{prop}\label{nochange}
 Let
$\lambda,\lambda^\prime\in\R$ such that $\lambda^\prime\leq\lambda$
and $[\lambda^\prime,\lambda]\cap\mathcal{D}=\emptyset$.  The
kernels, and cokernels, of \eq{ch6s3eq1} when $\mu=\lambda$ and
$\mu=\lambda^\prime$ are equal.
\end{prop}

\begin{proof}
Denote the dimensions of the kernel and cokernel of \eq{ch6s3eq1},
for $\mu\notin\mathcal{D}$, by $k(\mu)$ and $c(\mu)$ respectively.
Since $[\lambda^\prime,\lambda]\cap\mathcal{D}=\emptyset$,
$k(\lambda)-c(\lambda)=k(\lambda^\prime)-c(\lambda^\prime)$ and
hence
\begin{equation}\label{ch6s3subs2eq1a}
 k(\lambda)-k(\lambda^\prime)=c(\lambda)-c(\lambda^\prime).\end{equation}

We know that $k(\lambda)\leq k(\lambda^\prime)$ because
$L^p_{k+1,\,\lambda}\hookrightarrow L^p_{k+1,\,\lambda^\prime}$ by
Theorem \ref{ch6s2thm1} as $\lambda\geq\lambda^\prime$. Similarly,
since $c(\mu)$ is equal to the dimension of the kernel of the formal
adjoint operator acting on a Sobolev space with weight $-3-\mu$,
$c(\lambda)\geq c(\lambda^\prime)$.  Noting that the right-hand side
of \eq{ch6s3subs2eq1a} is non-negative and the left-hand side is
less than or equal to zero, we conclude that both must be zero. The
result follows from the fact that the kernel of $d+d^*$ in
$L^p_{k+1,\,\lambda}$ is contained in the kernel of $d+d^*$ in
$L^p_{k+1,\,\lambda^\prime}$, and vice versa for the cokernels.
\end{proof}

We can now go further and give a more explicit description of the
quantity $\d(\mu)$ in Definition \ref{ch6s3dfn1}.

\begin{prop}\label{ch6s3prop3}
Using the notation of Proposition \ref{ch6s3prop1} and Definition
\ref{ch6s3dfn1},
$\text{\emph{d}}(\mu)=\sum_{i=1}^s\text{\emph{dim}}\,D(\mu,i)$ for
$\mu\in\mathcal{D}$.
\end{prop}

\noindent This is an analogue of \cite[Proposition 5.4]{Lotay},
which is for the AC scenario, and can be proved in exactly the same
manner.

\section{The deformation problems}\label{ch8s2}

We have a common notation for the next three sections. Let $N$ be a
CS coassociative 4-fold at $z_1,\ldots,z_s$ with rate $\lambda$ in a
$\G2$ manifold $(M,\varphi,g)$. Suppose
$\lambda\in(1,2)\setminus\mathcal{D}$, where $\mathcal{D}$ is
defined in Proposition \ref{ch6s3prop1}, and the cone at $z_i$ is
$C_i$ with link $\Sigma_i$. We shall then use the notation of
Definitions \ref{ch8s1dfn2} and \ref{ch8s1dfn3}. In particular, we
let $\{\chi_i:B(0;\eta)\rightarrow V_i\,:\,i=1,\ldots,s\}$, with
$d\chi_i|_0=\zeta_i$ for all $i$, be the $\G2$ coordinate system
near $z_1,\ldots,z_s$ used to define $N$ and let $\hat{C}_i$ be the
tangent cone at $z_i$.  Recalling that $N$ is a CS manifold, in the
sense of Definition \ref{ch6s1dfn3}, we have a radius function
$\rho$ on $\Nhat$ as in Definition \ref{ch6s1dfn4}.

We consider deformations of $N$ which are CS coassociative 4-folds
at $s$ points with rate $\lambda$ in $(M,\varphi,g)$ with the same
cones at the singularities as $N$, but the singularities need not be
at the same points, nor have identical tangent cone. We also,
eventually, consider deforming the $\text{G}_2$ structure on the
ambient 7-manifold $M$.

\subsection{Problem 1: fixed singularities and $\G2$
structure}\label{ch8s2subs1}

The first deformation problem we consider is where the deformations
of $N$ have identical singular points to $N$ with the same rate,
cones and tangent cones, and the $\text{G}_2$ structure of $M$ is
fixed.

\begin{dfn}\label{ch8s2subs1dfn1}
The \emph{moduli space of deformations} $\mathcal{M}_1(N,\lambda)$
for Problem 1 is the set of $N^\prime$ in $(M,\varphi,g)$ which are
CS coassociative 4-folds  at $z_1,\ldots,z_s$ with rate $\lambda$,
having
cone $C_i$ and tangent cone
 $\hat{C}_i$ at $z_i$ for all $i$, such that
there exists a homeomorphism $h:N\rightarrow N^\prime$, isotopic to
the identity, with $h(z_i)=z_i$ for
$i=1,\ldots,s$ and such that 
$h|_{\Nhat}:\Nhat\rightarrow N^\prime\setminus\{z_1,\ldots,z_s\}$ is
a diffeomorphism.
\end{dfn}


We begin our formulation of a local description of
$\mathcal{M}_1(N,\lambda)$ with a result which is immediate from the
proof of \cite[Chapter IV, Theorem 9]{Lang} since $M$ is a
Riemannian manifold.
\begin{thm}\label{ch7s1thm1} Let $P$ be a closed embedded submanifold
of $M$.  There exist an open subset $V$ of the normal bundle
$\nu(P)$ of $P$ in $M$, containing the zero section, and an open set
$S$ in $M$ containing $P$, such that the exponential map
$\exp|_V:V\rightarrow S$ is a diffeomorphism.
\end{thm}

\begin{note} The proof of this result relies entirely on the observation
that $\exp|_{\nu(P)}$ is a local isomorphism upon the zero
section.\end{note}

This information provides us with a useful corollary.

\begin{cor}\label{ch8s2subs1cor1}
For $i=1,\ldots,s$ choose
$\Phi_i:(0,\epsilon)\times\Sigma_i\rightarrow
B(0;\eta)\subseteq\R^7$ uniquely by imposing the condition that
$$\Phi_i(r_i,\sigma_i)-\iota_i(r_i,\sigma_i)\in
(T_{r_i\sigma_i}\iota_i(C_i))^\perp$$ for all
$(r_i,\sigma_i)\in(0,\epsilon)\times\Sigma_i$, which can be achieved
by making $\epsilon$ smaller and $K$ larger if necessary.  Let
$P_i=\iota_i((0,\epsilon)\times\Sigma_i)$, $
Q_i=\Phi_i((0,\epsilon)\times\Sigma_i)$ and define
$n_i: \nu(P_i)\rightarrow\R^7$ by
$n_i(r_i\sigma_i,v)=v+\Phi_i(r_i,\sigma_i)$. 
For all $i$, there exist an open subset $\hat{V}_i$ of $\nu(P_i)$ in
$\R^7$, containing the zero section, and an open set $\hat{S}_i$ in
$B(0;\eta)\subseteq\R^7$ containing $Q_i$ such that
$n_i|_{\hat{V}_i}:\hat{V}_i\rightarrow\hat{S}_i$ is a
diffeomorphism.  Moreover, $\hat{V}_i$ and $\hat{S}_i$ can be chosen
to grow like $r_i$ on $(0,\epsilon)\times\Sigma_i$, for all $i$, and
such that $P_i\subseteq\hat{S}_i$.
\end{cor}

\begin{proof}
Note that $n_i$ takes the zero section of $\nu(P_i)$ to $Q_i$. By
the definition of $\Phi_i$, we see that $n_i$ is a local isomorphism
upon the zero section. Thus, the proof of Theorem \ref{ch7s1thm1}
gives open sets $\hat{V}_i$ and $\hat{S}_i$ such that
$n_i|_{\hat{V}_i}:\hat{V}_i\rightarrow\hat{S}_i$ is a
diffeomorphism. We can ensure that $\hat{S}_i$ lies in $B(0;\eta)$
by making $\hat{V}_i$ smaller if necessary.

Furthermore, since $\Phi_i-\iota_i$ is orthogonal to
$(0,\epsilon)\times\Sigma_i$, it can be identified with a small
section of the normal bundle and hence $P_i$ lies in $\hat{S}_i$ as
long as $\hat{S}_i$ grows at $O(r_i)$ as $r_i\rightarrow 0$.  As we
can form $\hat{S}_i$ and $\hat{V}_i$ in a translation equivariant
way because we are working on a portion of the cone $C_i$, we can
construct our sets with this decay rate as $r_i\rightarrow0$  and
such that they do not collapse as $r_i\rightarrow\epsilon$.

\end{proof}

Corollary \ref{ch8s2subs1cor1} helps us in establishing the next
proposition.

\begin{prop}\label{ch8s2subs1prop1}
There exist an open set $\hat{U}\subseteq\Lambda^2_+T^*\hat{N}$
containing the zero section, an open set $\hat{T}\subseteq M$
containing $\hat{N}$ and a diffeomorphism
$\delta:\hat{U}\rightarrow\hat{T}$ which takes the zero section to
$\hat{N}$. Moreover, $\hat{U}$ and $\hat{T}$ can be chosen to grow
with order $O(\rho)$ as $\rho\rightarrow 0$ and 
$\delta$ is compatible with the identifications
$U_i\setminus\{z_i\}\cong (0,\epsilon)\times\Sigma_i$ for all $i$
and the isomorphism
$\jmath:\nu(\Nhat)\rightarrow\Lambda^2_+T^*\Nhat$ given in
Proposition \ref{ch2s3subs2thm1}.
\end{prop}

\begin{proof} Use the notation of
Corollary \ref{ch8s2subs1cor1} and define
$\hat{T}_i=\chi_i(\hat{S}_i)$. Then $\hat{T}_i$ is an open set in
$M$ such that $U_i\setminus\{z_i\}\subseteq\hat{T}_i\subseteq V_i$,
since $\chi_i(Q_i)=U_i\setminus\{z_i\}$, and which grows with order
$O(\rho)$ as $\rho\rightarrow 0$.

Consider the bundle
$(\Lambda^2_+)_{\chi_i^*(g)}T^*((0,\epsilon)\times\Sigma_i)$, where
the notation $(\Lambda^2_+)_h$ indicates that the Hodge star is
calculated using the metric $h$ and we consider
$(0,\epsilon)\times\Sigma_i\cong P_i\subseteq\R^7$. Then
\begin{align*}
\jmath_i:\nu(P_i)&\longrightarrow(\Lambda^2_+)_{\chi_i^*(g)}T^*P_i\\
v|_{r_i\sigma_i}&\longmapsto\left(v|_{r_i\sigma_i}\cdot
 \chi_i^*(\varphi)|_{\Phi_i(r_i,\sigma_i)}\right)|_{T_{r_i\sigma_i}P_i}\end{align*}
is an isomorphism because $U_i\setminus\{z_i\}$ is coassociative and
thus $P_i$ is, with respect to the metric $\chi_i^*(g)$ and 3-form
$\chi_i^*(\varphi)$, and hence we may apply Proposition
\ref{ch2s3subs2thm1}. Note also that
$$\Psi_i^*:(\Lambda^2_+)_gT^*(U_i\setminus\{z_i\})\longrightarrow
(\Lambda^2_+)_{\chi_i^*(g)}T^*((0,\epsilon)\times\Sigma_i)$$ is
clearly a diffeomorphism. Therefore, let
$\hat{U}_i\subseteq(\Lambda^2_+)_gT^*(U_i\setminus\{z_i\})$ be such
that $\Psi_i^*(\hat{U}_i)=\jmath_i(\hat{V}_i)$.  Note, by
construction, that $\hat{U}_i$ grows with order $O(\rho)$ as
$\rho\rightarrow0$.

Define a diffeomorphism $\delta_i:\hat{U}_i\rightarrow\hat{T}_i$
such that the following diagram commutes:
\begin{equation}\label{ch8s2subs1eq1}
\begin{gathered}
\xymatrix{ 
\hat{U}_i\ar[rr]^{\Psi_i^*}\ar[dd]_{\delta_i} & &
\jmath_i(\hat{V}_i) \ar[d]^{\jmath_i^{-1}}
\\
&&\hat{V}_i\ar[d]^{n_i} 
\\
 \hat{T}_i & &\hat{S}_i.\ar[ll]_{\chi_i}
}\end{gathered}
\end{equation}

\noindent Interpolating smoothly over $K$, we extend
$\bigcup_{i=1}^s\hat{U}_i$ and $\bigcup_{i=1}^s\hat{T}_i$ to
$\hat{U}$ and $\hat{T}$ as required and extend the diffeomorphisms
$\delta_i$ smoothly to a diffeomorphism
$\delta:\Uhat\rightarrow\hat{T}$ such that $\delta$ acts as the
identity on $\Nhat$, which is identified with the zero section in
$\Lambda^2_+T^*\Nhat$.

Note that we have a splitting $T\Uhat|_{(x,0)}=T_x\Nhat\oplus
\Lambda^2_+T^*_x\Nhat$ for all $x\in\Nhat$.  Thus we can consider
$d\delta$ at $\Nhat$ as a map from $T\Nhat\oplus\Lambda^2_+T^*\Nhat$
to $T\Nhat\oplus\nu(\Nhat)\cong TM|_{\Nhat}$. Hence, we require in
our extension of $\delta$ from $\delta_i$ to ensure that, in matrix
notation,
\begin{equation}\label{ch8s2subs1eq2}d\delta|_{\Nhat}=\left(\begin{array}{cc}I& A\\
0 & \jmath^{-1}
\end{array}\right),\end{equation}
where $I$ is the identity and $A$ is arbitrary.  This can be
achieved because of the definition of $\delta_i$.

  The compatibility of $\delta$
with $\jmath\,$ and $\Psi_i$ for all $i$, mentioned in the statement
of the proposition, is given by \eq{ch8s2subs1eq1} and the behaviour
of $d\delta|_{\Nhat}$ stipulated in \eq{ch8s2subs1eq2}.

\end{proof}

We now define our deformation map for Problem 1.  Let
$C^k_{\text{loc}}(\hat{U})=\{\alpha\in
C^k_{\text{loc}}(\Lambda^2_+T^*\Nhat)\,:\,\alpha\in\Uhat\}$, where
$\Uhat$ is given in Proposition \ref{ch8s2subs1prop1}, and adopt
similar notation to define subsets of the spaces of forms described
in $\S$\ref{Banach}.

\begin{dfn}\label{ch8s2subs1dfn2}
Use the notation of Proposition \ref{ch8s2subs1prop1}.
 Let $\Gamma_\alpha$ be the graph of $\alpha\in C^1_{\text{loc}}(\hat{U})$ and
let $\pi_\alpha:\hat{N}\rightarrow\Gamma_\alpha$ be given by
$\pi_\alpha(x)=(x,\alpha(x))$.  Let $f_\alpha=\delta\circ\pi_\alpha$
and let $\hat{N}_{\alpha}=f_\alpha(\hat{N})\subseteq \hat{T}$.
Define a map $F_1$ from $C^1_{\text{loc}}(\hat{U})$ to
$C^0_{\text{loc}}(\Lambda^3T^*\hat{N})$ by:
$$F_1(\alpha)=f_\alpha^*\left(\varphi|_{\hat{N}_\alpha}\right).$$
By \cite[p. 731]{McLean}, which we are allowed to use by our choice
of $\delta$, the linearisation of $F_1$ at $0$ is
$$dF_1|_{0}(\alpha)=L_1(\alpha)=d\alpha$$
for all $\alpha\in C^1_{\text{loc}}(\Lambda^2_+T^*\Nhat)$.
\end{dfn}

\begin{remark}
The operator $L_1$ is \emph{not} elliptic.
\end{remark}

\noindent By Proposition \ref{ch2s3subs2prop2}, $\text{Ker}\,F_1$ is
the set of $\alpha\in C^1_{\text{loc}}(\hat{U})$ such that
$\hat{N}_\alpha$ is coassociative.

However, we want CS coassociative deformations with singularities at
the same points with the same tangent cones. Suppose $\alpha\in
C^1_{\text{loc}}(\Uhat)$ and
$N_{\alpha}=\Nhat_{\alpha}\cup\{z_1,\ldots,z_s\}$ is such a
deformation.  Then there exist smooth maps
$(\Phi_{\alpha})_i:(0,\epsilon)\times\Sigma_i \rightarrow B(0;\eta)$
satisfying \eq{ch8s1eq1} such that
$(\Psi_\alpha)_i=\chi_i\circ(\Phi_\alpha)_i$ is a diffeomorphism
onto an open subset of ${\Nhat_\alpha}$ for all $i$ as in Definition
\ref{ch8s1dfn2}.  Note that we are free to use $\chi_i$ because the
tangent cones at the singularities of $N_\alpha$ must be the same as
for $N$, so any $\text{G}_2$ coordinate system near the
singularities used to define $N_\alpha$ must be equivalent to the
one given by $\chi_i$ for $i=1,\ldots,s$. Choose $(\Phi_{\alpha})_i$
uniquely such that
$$(\Phi_{\alpha})_i(r_i,\sigma_i)-\iota_i(r_i,\sigma_i)\in
(T_{r_i\sigma_i}\iota_i(C_i))^\perp$$ for all
$(r_i,\sigma_i)\in(0,\epsilon)\times\Sigma_i$.

Use the notation of Corollary \ref{ch8s2subs1cor1} and the proof of
Proposition \ref{ch8s2subs1prop1}. Since
$$\Phi_i(r_i,\sigma_i)-\iota_i(r_i,\sigma_i)\in
(T_{r_i\sigma_i}P_i)^{\perp}\cong \nu_{r_i\sigma_i}(P_i),$$
$\Phi_i-\iota_i$ can be identified using $\jmath_i$ with the graph
of $\beta_i\in
(\Lambda_+^2)_{\chi_i^*(g)}T^*((0,\epsilon)\times\Sigma_i)$. Thus,
\begin{equation}\label{ch7s1eq13}
|\nabla_i^j\beta_i|=O(r_i^{\lambda-j})\qquad\text{for $j\in\N$ as
$r_i\rightarrow 0$}\end{equation} by \eq{ch8s1eq1} and therefore
$\beta_i\in
C_{\lambda}^{\infty}((\Lambda_+^2)_{\chi_i^*(g)}T^*((0,\epsilon)\times\Sigma_i))$.

We may similarly deduce, by the definition of $\delta$, $\Phi_i$ and
$(\Phi_{\alpha})_i$, that
$(\Phi_{\alpha})_i-\iota_i=((\Phi_{\alpha})_i-\Phi_i)+(\Phi_i-\iota_i)$
corresponds to the graph of $\Psi_i^*(\alpha)+\beta_i$ on
$(0,\epsilon)\times\Sigma_i$, recalling that
$$\Psi_i^*:\Lambda_+^2T^*(U_i\setminus\{z_i\})\rightarrow
(\Lambda_+^2)_{\chi_i^*(g)}T^*((0,\epsilon)\times\Sigma_i)$$ is a
diffeomorphism for all $i$.
Since $N_\alpha$ has the same types of singularities as $N$, both
$\beta_i$ and $\Psi_i^*(\alpha)+\beta_i$ lie in
$C_{\lambda}^{\infty}((\Lambda_+^2)_{\chi_i^*(g)}T^*((0,\epsilon)\times\Sigma_i))$
for each $i$.  Thus $\alpha$ must lie in
$C_{\lambda}^{\infty}(\Lambda_+^2T^*\Nhat)$.

We conclude that $\Nhat_{\alpha}$ is a sufficiently nearby
deformation of $\Nhat$ with the same conical singularities if and
only if $\alpha\in C_{\lambda}^{\infty}(\Uhat)\subseteq
C_{\lambda}^{\infty}(\Lambda_+^2T^*\Nhat)$.  We state this as a
proposition.

\begin{prop}\label{ch8s2subs1prop2}
The moduli space of deformations for Problem 1 is locally
homeomorphic to $\text{\emph{Ker}}\,F_1=\{\alpha\in
C^{\infty}_{\lambda}(\hat{U})\,:\, F_1(\alpha)=0\}$.
\end{prop}

We define an associated map $G_1$ to $F_1$.

\begin{dfn}\label{ch8s2subs1dfn3} Define $G_1:C_{\text{loc}}^{1}(\Uhat)\times
C_{\text{loc}}^{1}(\Lambda^4T^*\Nhat)\rightarrow
C^{0}_{\text{loc}}(\Lambda^3T^*\Nhat)$ by:
\begin{equation*}
G_1(\alpha,\beta)=F_1(\alpha)+d^*\beta. \end{equation*} Then $G_1$
is a first order \emph{elliptic} operator at $(0,0)$ since
$$dG_1|_{(0,0)}=d+d^*:C^1_{\text{loc}}(\Lambda^2_+T^*\Nhat\oplus\Lambda^4T^*\Nhat)
\longrightarrow C^0_{\text{loc}}(\Lambda^3T^*\Nhat).$$
\end{dfn}

\vspace{-28pt}

\begin{note} If $G_1(\alpha,\beta)=0$ and $\beta\in
C^{\infty}_{\lambda}(\Lambda^4T^*\hat{N})$, $*\beta$ is a harmonic
function which decays with order $O(\rho^{\lambda})$ as
$\rho\rightarrow 0$. Since $\lambda>1$, $*\beta\rightarrow 0$ as
$\rho\rightarrow 0$ and hence, by the Maximum Principle for harmonic
functions, it must be $0$.
\end{note}

We therefore deduce the following.

\begin{prop}\label{ch8s2subs1prop3}
$\text{\emph{Ker}}\,F_1\cong\{(\alpha,\beta)\in C^{\infty}_{\lambda}
(\hat{U})\times C^{\infty}_{\lambda}(\Lambda^4T^*\hat{N})\,:\,
G_1(\alpha,\beta)=0\}$.
\end{prop}

We conclude this subsection by stating and proving two results on
\emph{regularity} which are analogous to \cite[Proposition
4.3]{Lotay} and the argument in \cite[p. 22-24]{Lotay} respectively.

\begin{prop}\label{ch7s1prop3}
The map $F_1$ given in Definition \ref{ch8s2subs1dfn2} can be
written as
\begin{equation}\label{ch7s1eq3}
F_1(\alpha)(x)=d\alpha(x)+P_{F_1}(x,\alpha(x),\nabla\alpha(x))
\end{equation}
for $x\in \hat{N}$, where $P_{F_1}:\{(x,y,z)\,:\,(x,y)\in
\hat{U},\,z\in
T_x^*\hat{N}\otimes\Lambda_+^2T_x^*\hat{N})\}\rightarrow\Lambda^3T^*\hat{N}$
is a smooth map such that $P_{F_1}(x,y,z)\in\Lambda^3T^*_x\hat{N}$.
For $\alpha\in C_{\lambda}^{\infty}(\hat{U})$ with
$\|\alpha\|_{C^1_1}$ sufficiently small, denoting
$P_{F_1}(x,\alpha(x),\nabla\alpha(x))$ by $P_{F_1}(\alpha)(x)$,
$P_{F_1}(\alpha)\in
C_{2\lambda-2}^{\infty}(\Lambda^3T^*\hat{N})\subseteq
C_{\lambda-1}^{\infty}(\Lambda^3T^*\hat{N})$, as $\lambda>1$.
Moreover, for each $k\in\N$, if $\alpha\in C^{k+1}_\lambda(\hat{U})$
and $\|\alpha\|_{C^1_1}$ is sufficiently small, $P_{F_1}(\alpha)\in
C^k_{2\lambda-2}(\Lambda^3T^*\hat{N})$ and there exists a constant
$c_k>0$ such that
$$\|P_{F_1}(\alpha)\|_{C_{2\lambda-2}^k}\leq
c_k\|\alpha\|_{C^{k+1}_{\lambda}}^2\,.
$$
\end{prop}
\begin{proof} Firstly, by the definition of $F_1$, $F_1(\alpha)(x)$
relates to the tangent space to $\Gamma_{\alpha}$ at
$\pi_{\alpha}(x)$.  Note that $T_{\pi_\alpha(x)}\Gamma_{\alpha}$
depends on both $\alpha(x)$ and $\nabla\alpha(x)$ and hence so must
$F_1(\alpha)(x)$.  We may then define $P_{F_1}$ by \eq{ch7s1eq3}
such that it is a smooth function of its arguments as claimed.

We argued above that we may identify $\Phi_i-\iota_i$ on
$(0,\epsilon)\times\Sigma_i$ with $$\beta_i\in
C^{\infty}_{\lambda}((\Lambda_+^2)_{\chi_i^*(g)}T^*((0,\epsilon)\times\Sigma_i))$$
for $i=1,\ldots,s$.  Recall that
$$\Psi_i^*:\Lambda_+^2T^*(U_i\setminus\{z_i\})\rightarrow
(\Lambda_+^2)_{\chi_i^*(g)}T^*((0,\epsilon)\times\Sigma_i)$$ is a
diffeomorphism.  Let $k\in\N$, $\alpha\in C^{k+1}_\lambda(\hat{U})$,
$\alpha_i=\alpha|_{U_i\setminus\{z_i\}}$ and
$\gamma_i=\Psi_i^*(\alpha_i)$.

For each $i$, define a function $F_{C_i}(\gamma_i+\beta_i)$ on
$(0,\epsilon)\times\Sigma_i$ by
\begin{equation}\label{ch7s1eq4}
F_{C_i}(\gamma_i+\beta_i)(r_i,\sigma_i)=F_1(\alpha_i)(\Psi_i(r_i,\sigma_i)).
\end{equation}
Define a smooth function $P_{C_i}$ by an equation analogous to
\eq{ch7s1eq3}:
\begin{align}
F_{C_i}(\gamma_i+\beta_i)(r_i,\sigma_i)&=
d(\gamma_i+\beta_i)(r_i,\sigma_i)\nonumber\\&+P_{C_i}(
(r_i,\sigma_i),(\gamma_i+\beta_i)(r_i,\sigma_i),\nabla(\gamma_i+\beta_i)(r_i,\sigma_i)).\label{ch7s1eq5}
\end{align}
We notice that $F_{C_i}$ and $P_{C_i}$ are only dependent on the
cone $C_i$ and, rather trivially, on $\epsilon$.  Therefore, because
of this fact and our choice of $\delta$ in Proposition
\ref{ch8s2subs1prop1}, these functions have scale equivariance
properties.  We may therefore derive equations and inequalities on
$\{\epsilon\}\times\Sigma_i$ and deduce the result on all of
$(0,\epsilon)\times\Sigma_i$ by introducing an appropriate scaling
factor of $r$.

Now, since $\alpha=0$ corresponds to our coassociative 4-fold
$\hat{N}$, $F_1(0)=0$. So, by \eq{ch7s1eq4},
\begin{equation}\label{ch7s1eq14}
F_{C_i}(\beta_i)=d\beta_i+P_{C_i}(\beta_i)=0, \end{equation}
adopting similar notation for $P_{C_i}(\beta_i)$ as for
$P_{F_1}(\alpha_i)$. Using \eq{ch7s1eq3}-\eq{ch7s1eq14}, we deduce
that
\begin{align}
P_{F_1}(\alpha_i)&
=d\beta_i+P_{C_i}(\gamma_i+\beta_i) 
=d\beta_i+P_{C_i}(\gamma_i+\beta_i)-(d\beta_i+P_{C_i}(\beta_i)) \nonumber\\
\label{ch7s1eq6}
&=P_{C_i}(\gamma_i+\beta_i)-P_{C_i}(\beta_i).
\end{align}
We then calculate
\begin{align}
P_{C_i}(\gamma_i+\beta_i)-P_{C_i}(\beta_i)&=\int_0^1\frac{d}{dt}\,P_{C_i}(t\gamma_i+\beta_i)
\,dt
\nonumber\\[4pt]
\label{ch7s1eq7}
 & =\int_0^1\!\gamma_i\cdot\frac{\partial P_{C_i}}{\partial
y}(t\gamma_i+\beta_i)+\nabla\gamma_i\cdot\frac{\partial
P_{C_i}}{\partial z}(t\gamma_i+\beta_i)\,dt,
\end{align}
recalling that $P_{C_i}$ is a function of three variables $x$, $y$
and $z$.
 Using Taylor's Theorem,
\begin{equation}\label{ch7s1eq8}
P_{C_i}(\gamma_i+\beta_i)=P_{C_i}(\beta_i)+\gamma_i\cdot\frac{\partial
P_{C_i}}{\partial y}(\beta_i)+ \nabla\gamma_i\cdot\frac{\partial
P_{C_i}}{\partial z}(\beta_i)
+O(r^{-2}|\gamma_i|^2+|\nabla\gamma_i|^2)\end{equation} when
$|\gamma_i|$ and $|\nabla\gamma_i|$ are small. Since
$dF_1|_0(\alpha_i)=d\alpha_i$,
$dF_{C_i}|_{\beta_i}(\gamma_i+\beta_i)=d\gamma_i$ and hence
$dP_{C_i}|_{\beta_i}=0$.  Thus, the first derivatives of $P_{C_i}$
with respect to $y$ and $z$ must vanish at $\beta_i$ by
\eq{ch7s1eq8}.  Therefore, given small $\nu>0$ there exists a
constant $A_0>0$ such that
\begin{equation}\label{ch7s1eq9}
\begin{split}\left|\frac{\partial P_{C_i}}{\partial
y}(t\gamma_i+\beta_i)\right| &\leq A_0(r^{-2}|\gamma_i|
+r^{-1}|\nabla\gamma_i|); \;\text{and}\\
\left|\frac{\partial P_{C_i}}{\partial
z}(t\gamma_i+\beta_i)\right|&\leq A_0(r^{-1}|\gamma_i|
+|\nabla\gamma_i|)\end{split}\end{equation} for $t\in [0,1]$
whenever
\begin{equation}\label{ch7s1eq10}
r^{-1}|\gamma_i|,\,r^{-1}|\beta_i|,\,|\nabla\gamma_i|\;\text{and}\;|\nabla\beta_i|\leq\nu.
\end{equation}
By \eq{ch7s1eq13}, $r^{-1}|\beta_i|$ and $|\nabla\beta_i|$ tend to
zero as $r\rightarrow0$. We can thus ensure that \eq{ch7s1eq10} is
satisfied by the $\beta_i$ components by making $\epsilon$ smaller.
Hence, \eq{ch7s1eq10} holds if
$\|\gamma_i\|_{C^1_1}\leq\nu$.
 Therefore, putting estimates
\eq{ch7s1eq9} in \eq{ch7s1eq7} and using \eq{ch7s1eq6},
\begin{equation}\label{ch7s1eq11}
|P_{F_1}(\alpha_i)|=|P_{C_i}(\gamma_i+\beta_i)-P_{C_i}(\beta_i)|\leq
A_0(r^{-1}|\gamma_i|+|\nabla\gamma_i|)^2
\end{equation}
whenever $\|\gamma_i\|_{C^1_1}\leq\nu$.  As $r\rightarrow0$ the
terms in the bracket on the right-hand side of \eq{ch7s1eq11} are of
order $O(r^{\lambda-1})$ by \eq{ch7s1eq13}. Thus,
$|P_{F_1}(\alpha_i)|$ is of order $O(r^{2\lambda-2})$, hence
$O(r^{\lambda-1})$ since $\lambda>1$, as $r\rightarrow0$ for
$i=1,\ldots,s$. We deduce that $|P_{F_1}(\alpha)|$ is of order
$O(\rho^{2\lambda-2})$ as $\rho\rightarrow 0$ for all $\alpha\in
C^1_\lambda(\hat{U})$ with $\|\alpha\|_{C^1_1}$ sufficiently small.

Similar calculations give analogous results to \eq{ch7s1eq11} for
derivatives of $P_{F_1}$, but we shall explain the method by
considering the first derivative. From \eq{ch7s1eq7} we calculate
\begin{align*}
\nabla&(P_{C_i}(\gamma_i+\beta_i)-P_{C_i}(\beta_i))\\[4pt]=&\int_0^1
\nabla\!\left(\gamma_i\cdot\frac{\partial P_{C_i}}{\partial
y}(t\gamma_i+\beta_i)+\nabla\gamma_i\cdot\frac{\partial
P_{C_i}}{\partial z}(t\gamma_i+\beta_i)\right)dt \\[4pt]
=&\int_0^1\nabla\gamma_i\cdot\frac{\partial P_{C_i}}{\partial
y}+\gamma_i\cdot\!\left(\nabla(t\gamma_i+\beta_i)\cdot\frac{\partial^2
P_{C_i}}{\partial
y^2}+\nabla^2(t\gamma_i+\beta_i)\cdot\frac{\partial^2
P_{C_i}}{\partial y\partial z}\right) \\[4pt]
&+\nabla^2\gamma_i\cdot\frac{\partial P_{C_i}}{\partial
z}+\nabla\gamma_i\cdot\!\left(\nabla(t\gamma_i+\beta_i)\cdot\frac{\partial^2
P_{C_i}}{\partial z\partial
y}+\nabla^2(t\gamma_i+\beta_i)\cdot\frac{\partial^2
P_{C_i}}{\partial z^2}\right)dt.
\end{align*}
Whenever $\|\gamma_i\|_{C^1_1}\leq\nu$ there exists a constant
$A_1>0$ such that \eq{ch7s1eq9} holds with $A_0$ replaced by $A_1$
and, for $t\in[0,1]$,
$$\left|\frac{\partial^2
P_{C_i}}{\partial
y^2}(t\gamma_i+\beta_i)\right|,\,\left|\frac{\partial^2
P_{C_i}}{\partial y\partial
z}(t\gamma_i+\beta_i)\right|\;\text{and}\; \left|\frac{\partial^2
P_{C_i}}{\partial z^2}(t\gamma_i+\beta_i)\right|\leq A_1,$$ since
the second derivatives of $P_{C_i}$ are continuous functions defined
on the closed
bounded set given by $\|\gamma_i\|_{C^1_1}\leq\nu$.  
We deduce that
$$
\big|\nabla
\big(P_{F_1}(\alpha_i)\big)\big|=|\nabla(P_{C_i}(\gamma_i+\beta_i)-P_{C_i}(\beta_i))|\leq
A_1\left(\sum_{j=0}^2r^{j-2}|\nabla^j\gamma_i| \right)^{\!2}
$$
whenever $\|\gamma_i\|_{C^1_1}\leq\nu$.  Therefore $\big|\nabla
\big(P_{F_1}(\alpha_i)\big)\big|$ is of order $O(r^{2\lambda-3})$,
hence $O(r^{\lambda-2})$, as $r\rightarrow 0$.

In general we have the estimate
$$\big|\nabla^l\big(P_{F_1}(\alpha_i)\big)\big|\leq A_l
\left(\sum_{j=0}^{l+1}r^{j-(l+1)}|\nabla^j\gamma_i|\right)^{\!2}$$
for some $A_l>0$ whenever $\|\gamma_i\|_{C^1_1}\leq\nu$.  The result
follows.
\end{proof}

We now consider the regularity of solutions to the \emph{nonlinear}
elliptic equation $G_1(\alpha,\beta)=0$ near $(0,0)$.

\begin{prop}\label{ch8s2subs1prop4} Let $(\alpha,\beta)\in
L_{k+1,\,\lambda}^p(\Uhat)\times
L_{k+1,\,\lambda}^p(\Lambda^4T^*\Nhat)$ for some $p>4$ and $k\geq
2$. If $G_1(\alpha,\beta)=0$ and $\|\alpha\|_{C^1_1}$ is
sufficiently small, $(\alpha,\beta)\in
C_{\lambda}^{\infty}(\Uhat)\times
C_{\lambda}^{\infty}(\Lambda^4T^*\Nhat)$.
\end{prop}

\begin{proof}
Suppose that $(\alpha,\beta)\in L_{k+1,\,\lambda}^p(\Uhat)\times
L_{k+1,\,\lambda}^p(\Lambda^4T^*\Nhat)$ for some $p>4$ and $k\geq
2$. Then $\alpha$ and $\beta$ lie in $C^1_{\text{loc}}$ by Theorem
\ref{ch6s2thm1}, since $\frac{k}{4}>\frac{1}{p}$\,.

Suppose further that $G_1(\alpha,\beta)=0$ and that
$\|\alpha\|_{C^1_1}$ is sufficiently small. Since $F_1$ smoothly
depends on $\alpha$ and $\nabla\alpha$, $G_1$ is a smooth function
of $\alpha,\beta,\nabla\alpha$ and $\nabla\beta$. We apply
\cite[Theorem 6.8.1]{Morrey}, which is a general regularity result
for nonlinear elliptic equations, to conclude that $\alpha$ and
$\beta$ are smooth. However, we want more than this: the derivatives
of $\alpha$ and $\beta$ must decay at the required rates.

Recall the note after Definition \ref{ch8s2subs1dfn3} that
$G_1(\alpha,\beta)=0$ implies that $\beta=0$. Thus $\beta\in
C^{\infty}_\lambda(\Lambda^4T^*\Nhat)$ trivially.

For the following argument we find it useful to work with weighted
H\"older spaces.  By Theorem \ref{ch6s2thm1}, $\alpha\in
C^{k,\,a}_{\lambda}(\Uhat)$ with $a=1-4/p\in (0,1)$ since $p>4$. We
also know that $d^*(G_1(\alpha,\beta))=d^*(F_1(\alpha))=0$, which is
a nonlinear elliptic equation on $\alpha$.  Using the notation and
results of Proposition \ref{ch7s1prop3},
$d^*d\alpha+d^*(P_{F_1}(\alpha))=0$ and $d^*(P_{F_1}(\alpha))\in
C^{k-2,\,a}_{2\lambda-3}(\Lambda^2T^*\Nhat)$.  We see that
$$d^*(F_1(\alpha))(x)=R(x,\alpha(x),\nabla\alpha(x))\nabla^2\alpha(x)
+E(x,\alpha(x),\nabla\alpha(x)),$$ where
$R(x,\alpha(x),\nabla\alpha(x))$ and
$E(x,\alpha(x),\nabla\alpha(x))$ are smooth functions of their
arguments.  Define
$$S_{\alpha}(\gamma)(x)= R(x,\alpha(x),\nabla\alpha(x))
\nabla^2\gamma(x)$$ for $\gamma\in
C^2_{\text{loc}}(\Lambda^2_+T^*\Nhat)$. Then $S_{\alpha}$ is a
smooth, \emph{linear}, \emph{elliptic}, second-order operator, if
$\|\alpha\|_{C^1_1}$ is sufficiently small, whose coefficients
depend on $x$, $\alpha(x)$ and $\nabla\alpha(x)$. These coefficients
therefore lie in $C^{k-1,\,a}_\text{loc}$.  We also notice that
$$S_{\alpha}(\alpha)(x)=-E(x,\alpha(x),\nabla\alpha(x))\in
C^{k-2,\,a}_{2\lambda-3}(\Lambda^2T^*\Nhat)\subseteq
C^{k-2,\,a}_{\lambda-2}(\Lambda^2T^*\Nhat),$$ since $\lambda>1$.
However, $E(x,\alpha(x),\nabla\alpha(x))$ only depends on $\alpha$
and $\nabla\alpha$, and is at worst quadratic in these quantities by
Proposition \ref{ch7s1prop3}, so it must in fact lie in
$C^{k-1,\,a}_{\lambda-2}(\Lambda^2T^*\Nhat)$ since we are given
control on the decay of the first $k$ derivatives of $\alpha$ as
$\rho\rightarrow 0$.

The work in \cite[$\S$6.1.1]{Marshall} on asymptotically conical
manifolds gives regularity results for smooth linear elliptic
operators acting between weighted H\"older spaces.  These results
can easily be adapted to the CS scenario.  In particular, if
$\gamma\in C^2_{\lambda}(\Lambda_+^2T^*N)$ and
$S_{\alpha}(\gamma)\in C^{k-1,\,a}_{\lambda-2}(\Lambda^2T^*N)$, we
have that $\gamma\in C^{k+1,\,a}_{\lambda}(\Lambda_+^2T^*N)$.  Since
$k\geq 2$, $\alpha$ and $S_{\alpha}(\alpha)$ satisfy these
conditions by the discussion above. We deduce that $\alpha\in
C^{k+1,\,a}_{\lambda}(\Lambda^2_+T^*N)$ only knowing a priori that
$\alpha\in C^{k,\,a}_{\lambda}(\Lambda^2_+T^*N)$. We proceed by
induction to show that $\alpha\in
C^{l,\,a}_{\lambda}(\Lambda^2_+T^*N)$ for all $l\geq 2$.
\end{proof}

\subsection{Problem 2: moving singularities and fixed $\G2$
structure}\label{ch8s2subs2}

For this problem we again consider deformations of $N$ in
$(M,\varphi,g)$ which are CS coassociative 4-folds at $s$ points
with the same rate and cones at the singularities, but now we allow
the singular points and tangent cones at those points to differ from
those of $N$. However, we still assume that the $\text{G}_2$
structure on $M$ is fixed.

\begin{dfn}\label{ch8s2subs2dfn1}
The \emph{moduli space of deformations} $\mathcal{M}_2(N,\lambda)$
for Problem 2 is the set of $N^\prime$ in $(M,\varphi,g)$ which are
CS coassociative 4-folds at $z_1^\prime,\ldots,z_s^\prime$ with rate
$\lambda$, having cone $C_i$ and tangent cone
 $\hat{C}_i^\prime$ at $z_i^\prime$ for all $i$, such that
there exists a homeomorphism $h:N\rightarrow N^\prime$, isotopic to
the identity, with $h(z_i)=z_i^\prime$ for $i=1,\ldots,s$ and such
that $h|_{\Nhat}:\Nhat\rightarrow
N^\prime\setminus\{z_1^\prime,\ldots,z_s^\prime\}$ is a
diffeomorphism. 
\end{dfn}

Here it is more difficult to create a local description of the
moduli space which is compatible with the analytic framework in
which our study is made.  What one would consider more `intuitive'
approaches do not, as far as the author is aware, bear fruit.  We
therefore follow what is, at first sight, a slightly indirect route.

\medskip

For each $i=1,\ldots,s$ let $B_i$ be an open set in $M$ containing
$z_i$ such that $B_i\cap B_j=\emptyset$ for $i\neq j$. Let
$B=\prod_{i=1}^s B_i$.   For each
$\mathbf{z}^\prime=(z_1^\prime,\ldots,z_s^\prime)\in B$,
 we have a family $I(\mathbf{z}^\prime)$ of choices of
$s$-tuples
$\boldsymbol{\zeta}^\prime=(\zeta_1^\prime,\ldots,\zeta_s^\prime)$
of isomorphisms $\zeta_i^\prime:\R^7\rightarrow T_{z_i^\prime}M$
identifying $(\varphi_0,g_0)$ with
$(\varphi|_{T_{z_i^\prime}M},g|_{T_{z_i^\prime}M})$. Clearly, for
each $\mathbf{z}^\prime\in B$,
$I(\mathbf{z}^\prime)\cong\text{G}_2$.  We thus make the following
definition.

\begin{dfn}\label{ch8s2subs2dfn2}
The \emph{translation space} is
$$\mathcal{T}=\{(\mathbf{z}^\prime,\boldsymbol{\zeta}^\prime)\,:\,
\mathbf{z}^\prime\in B,\,\boldsymbol{\zeta}^\prime\in I
(\mathbf{z}^\prime)\}.
$$
It is a principal $\text{G}_2^s$ bundle over $B$ and hence is a
smooth manifold.

Let $\text{H}_i$ denote the Lie subgroup of $\text{G}_2$ preserving
$\iota_i(C_i)$ in $\R^7$ for $i=1,\ldots,s$ and let
$\text{H}=\prod_{i=1}^s\text{H}_i\subseteq\text{G}_2^s$.  Then
$\text{H}$ acts freely on $\mathcal{T}$ by
$$(\mathbf{z}^\prime,\boldsymbol{\zeta}^\prime)\longmapsto
(\mathbf{z}^\prime,(\zeta_1^\prime\circ
A_1^{-1},\ldots,\zeta_s^\prime\circ A_s^{-1})),$$ where
$(A_1,\ldots,A_s)\in\text{H}$. Thus there exists an $\text{H}$-orbit
through $(\mathbf{z},\boldsymbol{\zeta})$ in $\mathcal{T}$, where
$$\mathbf{z}=(z_1,\ldots,z_s)\quad\text{and}\quad\boldsymbol{\zeta}=
(\zeta_1,\ldots,\zeta_s).$$ Define $\hat{\mathcal{T}}$ to be a small
open ball in $\R^n$ containing $0$, where $n=
\text{dim}\,\mathcal{T}-\text{dim}\,\text{H}$, and let
$h_{\hat{\mathcal{T}}}:\hat{\mathcal{T}}\rightarrow\mathcal{T}$ be
an embedding with
$h_{\hat{\mathcal{T}}}(0)=(\mathbf{z},\boldsymbol{\zeta})$ such that
$h_{\hat{\mathcal{T}}}(\hat{\mathcal{T}})$ is transverse to the
$\text{H}$-orbit through $(\mathbf{z},\boldsymbol{\zeta})$. 
Write
$h_{\hat{\mathcal{T}}}(t)=(\mathbf{z}(t),\boldsymbol{\zeta}(t))$ for
$t\in\hat{\mathcal{T}}$, with $\mathbf{z}(0)=\mathbf{z}$ and
$\boldsymbol{\zeta}(0)=\boldsymbol{\zeta}$.
\end{dfn}

\begin{notes}\begin{itemize}\item[]\item[(a)] 
If $t,t^\prime\in\hat{\mathcal{T}}$, with $t\neq t^\prime$, are such
that $\mathbf{z}(t)=\mathbf{z}(t^\prime)$, the $s$-tuples of tangent
cones, $\{\hat{C}_1(t),\ldots,\hat{C}_s(t)\}$ and
$\{\hat{C}_1(t^\prime),\ldots,\hat{C}_s(t^\prime)\}$, are distinct.
\item[(b)] $\hat{\mathcal{T}}$ is an open ball in $\R^n\cong
T_0\hat{\mathcal{T}}$ and hence can be considered as an open subset
of $T_0\hat{\mathcal{T}}$.\end{itemize}\end{notes}

We use $\hat{\mathcal{T}}$ to extend $N$ to a family of nearby CS
4-folds and provide an analogue to Proposition \ref{ch8s2subs1prop1}
for Problem 2.  In defining $N$ we chose a $\text{G}_2$ coordinate
system $\{\chi_i:B(0;\eta)\rightarrow V_i:i=1,\ldots,s\}$ with
$d\chi_i|_0=\zeta_i$ for $i=1,\ldots,s$.  Extend this to a smooth
family of $\text{G}_2$ coordinate systems
$$\left\{\{\chi_i(t):B(0;\eta)\rightarrow V_i(t):i=1,\ldots,s\}:t\in
\hat{\mathcal{T}}\right\},$$ where $V_i(t)$ is an open set in $M$
containing $z_i(t)$, $\chi_i(t)(0)=z_i(t)$,
$d\chi_i(t)|_0=\zeta_i(t)$, $\chi_i(0)=\chi_i$ and $V_i(0)=V_i$ for
$i=1,\ldots,s$.

\begin{prop}\label{ch8s2subs2prop1} Use the notation of Proposition
\ref{ch8s2subs1prop1} and Definition \ref{ch8s2subs2dfn2}.
\begin{itemize}\item[(a)] There exists a family
$\mathcal{N}=\{N(t):t\in\hat{\mathcal{T}}\}$ of CS 4-folds in $M$,
with $N(0)=N$, such that $N(t)$ has singularities at
$z_1(t),\ldots,z_s(t)$ with rate $\lambda$, cones $C_1,\ldots,C_s$
and tangent cones $\hat{C}_1(t),\ldots,\hat{C}_s(t)$ defined by
$\hat{C}_i(t)=(\zeta_i(t)\circ\iota_i)(C_i)$.
\item[(b)] Let
$\hat{N}(t)=N(t)\setminus\{z_1(t),\ldots,z_s(t)\}$ and write
$$N(t)=K(t)\sqcup\bigsqcup_{i=1}^sU_i(t)$$
where $K(t)$ is compact and
$U_i(t)\setminus\{z_i(t)\}\cong(0,\epsilon)\times\Sigma_i$ for all
$i$, in the obvious way, ensuring that $K(0)=K$ and $U_i(0)=U_i$.
For $t\in\hat{\mathcal{T}}$, there exist open sets
$\hat{T}(t)\subseteq M$ containing $\hat{N}(t)$ and diffeomorphisms
$\delta(t):\hat{U}\rightarrow\hat{T}(t)$ taking the zero section to
$\hat{N}(t)$, varying smoothly in $t$, with $\hat{T}(0)=\hat{T}$ and
$\delta(0)=\delta$.  Moreover, $\hat{T}(t)$ can be chosen to grow
with order $O(\rho)$ as $\rho\rightarrow 0$ and $\delta(t)$ is
compatible with the identifications $U_i(t)\setminus\{z_i(t)\}\cong
(0,\epsilon)\times\Sigma_i$ for all $i$.
\end{itemize}
\end{prop}

\begin{remark} The family $\mathcal{N}$ does \emph{not} necessarily consist of
CS \emph{coassociative} 4-folds and $\delta(t)$ is not required to
be compatible with the isomorphism
$\nu(\Nhat)\cong\Lambda^2_+T^*\Nhat$ for $t\neq 0$.\end{remark}

\begin{proof}
Use the notation from the proof of Proposition
\ref{ch8s2subs1prop1}.  For $t\in\hat{\mathcal{T}}$, define
$\hat{T}_i(t)=\chi_i(t)(\hat{S}_i)$ and
$$U_i(t)=\Big(\chi_i(t)\circ\Phi_i((0,\epsilon)\times\Sigma_i)\Big)\cup\{z_i(t)\}$$
for $i=1,\ldots,s$.
 Then $\hat{T}_i(t)$ contains
$U_i(t)\setminus\{z_i(t)\}$.  Define a diffeomorphism $\delta_i(t)$
such that the following diagram commutes:
\begin{equation}\label{ch8s2subs2eq1}
\begin{gathered}\xymatrix{ 
\hat{U}_i\ar[rr]^{\Psi_i^*}\ar[dd]_{\delta_i(t)} & &
\jmath_i(\hat{V}_i) \ar[d]^{\jmath_i^{-1}}
\\
&&\hat{V}_i\ar[d]^{n_i} 
\\
 \hat{T}_i(t) & &\hat{S}_i.\ar[ll]_{\chi_i(t)}
}\end{gathered}\end{equation}

\noindent We then interpolate smoothly over $K$ to extend
$\bigcup_{i=1}^s\hat{T}_i(t)$ to $\hat{T}(t)$ and $\delta_i(t)$ to
$\delta(t)$ as required. Note by construction that $\hat{T}(t)$
grows with order $O(\rho)$ as $\rho\rightarrow0$.

Let $e(t)=\delta(t)|_{\Nhat}$ and define $\Nhat(t)=e(t)(\Nhat)$.
Then $e(t):\Nhat\rightarrow\Nhat(t)$ is a diffeomorphism for all
$t\in\hat{\mathcal{T}}$ and $e(0)$ is the identity.  Let
$N(t)=\Nhat(t)\cup\{z_1(t),\ldots,z_s(t)\}$.  We then have a family
$\mathcal{N}=\{N(t):t\in\hat{\mathcal{T}}\}$ as claimed. 
Note that $K(t)=e(t)(K)$.


By the construction of $\delta(t)$ and the family $\mathcal{N}$, it
is clear that the proposition is proved, where the compatibility
conditions on $\delta(t)$ are given by \eq{ch8s2subs2eq1}.
\end{proof}

 The next definition is analogous to Definition
\ref{ch8s2subs1dfn2}.

\begin{dfn}\label{ch8s2subs2dfn3}
Use the notation of Proposition \ref{ch8s2subs2prop1}.  Let
$\Gamma_\alpha$ be the graph of $\alpha\in C^1_{\text{loc}}(\Uhat)$
and let $\pi_\alpha:\Nhat\rightarrow\Gamma_\alpha$ be given by
$\pi_\alpha(x)=(x,\alpha(x))$.  For $t\in\hat{\mathcal{T}}$, let
$f_\alpha(t)=\delta(t)\circ\pi_\alpha$ and let
$\Nhat_\alpha(t)=f_\alpha(t)(\Nhat)$.  Define $F_2$ from
$C^1_{\text{loc}}(\Uhat)\times\hat{\mathcal{T}}$ to
$C^0_{\text{loc}}(\Lambda^3T^*\Nhat)$ by:
$$F_2(\alpha,t)=f_\alpha(t)^*\left(\varphi|_{\Nhat_\alpha(t)}\right).$$
\noindent The linearisation of $F_2$ at $(0,0)$ acts as
$$dF_2|_{(0,0)}:(\alpha,t)\longmapsto d\alpha+L_2(t),$$
where $\alpha\in C^1_{\text{loc}}(\Lambda^2_+T^*\Nhat)$, $t\in
T_0\hat{\mathcal{T}}$ and $L_2$ is a linear map into the space of
smooth exact 3-forms on $\Nhat$ since $\varphi$ is exact near
$\Nhat$.
\end{dfn}

\begin{remark} By construction $F_2(\alpha,0)=F_1(\alpha)$ as given in Definition
\ref{ch8s2subs1dfn2}. \end{remark}

Clearly, $\text{Ker}\,F_2$ is the set of $\alpha\in
C^1_{\text{loc}}(\Uhat)$ and $t\in\hat{\mathcal{T}}$ such that
$\Nhat_\alpha(t)$ is coassociative.  However, we have not yet
encoded the information that $N_\alpha(t)$ is CS with \emph{rate}
$\lambda$. This is the subject of the next proposition.

\begin{prop}\label{ch8s2subs2prop2}
The moduli space of deformations for Problem 2 is locally
homeomorphic to $\text{\emph{Ker}}\,F_2=\{(\alpha,t)\in
C^{\infty}_{\lambda}(\Uhat)\times\hat{\mathcal{T}}\,:\,
F_2(\alpha,t)=0\}$.
\end{prop}

\begin{proof} For each $t\in\hat{\mathcal{T}}$, we are in the situation of Problem
1 in the sense that we want coassociative deformations
$\Nhat_\alpha(t)$ of $\Nhat(t)$, defined by a self-dual 2-form
$\alpha$, which have the same singular points, cones and tangent
cones as $\Nhat(t)$.  It is thus clear that $\alpha\in
C^\infty_\lambda(\Uhat)$ by Proposition \ref{ch8s2subs1prop2}.
\end{proof}

We now introduce an associated map $G_2$ to $F_2$.

\begin{dfn}\label{ch8s2subs2dfn4} Define $G_2:C_{\text{loc}}^{1}(\Uhat)\times
C_{\text{loc}}^{1}(\Lambda^4T^*\Nhat)\times\hat{\mathcal{T}}\rightarrow
C^{0}_{\text{loc}}(\Lambda^3T^*\Nhat)$ by:
\begin{equation*}
G_2(\alpha,\beta,t)=F_2(\alpha,t)+d^*\beta. \end{equation*} Then
$dG_2|_{(0,0,0)}:(\alpha,\beta,t)\longmapsto
d\alpha+d^*\beta+L_2(t)$, in the notation of Definition
\ref{ch8s2subs2dfn3}.
\end{dfn}

\noindent We then have an analogous result to Proposition
\ref{ch8s2subs1prop3}, which follows in exactly the same fashion
because $F_2(\alpha,t)$ is exact.

\begin{prop}\label{ch8s2subs2prop3}
$$\text{\emph{Ker}}\,F_2\cong\{(\alpha,\beta,t)\in
C^\infty_\lambda(\Uhat)\times
C^{\infty}_{\lambda}(\Lambda^4T^*\Nhat)\times\hat{\mathcal{T}}\,:\,G_2(\alpha,\beta,t)=0\}.$$
\end{prop}

The next result studies the regularity of the kernel of $G_2$ near
$(0,0,0)$ and is the analogue of Proposition \ref{ch8s2subs1prop4}.

\begin{prop}\label{ch8s2subs2prop4}
Let $(\alpha,\beta,t)\in L^p_{k+1,\,\lambda}(\Uhat)\times
L^p_{k+1,\,\lambda}(\Lambda^4T^*\Nhat)\times\hat{\mathcal{T}}$,
where\linebreak $p>4$ and $k\geq 2$. If $G_2(\alpha,\beta,t)=0$ and
$\|\alpha\|_{C^1_1}$ and $t$ are sufficiently small,
$(\alpha,\beta)\in C^{\infty}_\lambda(\Uhat)\times
C^{\infty}_\lambda(\Lambda^4T^*\Nhat)$.
\end{prop}

\begin{proof}
Note that $dG_2(\alpha,\beta,t)=\Delta\beta=0$ implies that
$\beta=0$ by the Maximum Principle for harmonic functions and
$d^*G_2(\alpha,\beta,t)=d^*F_2(\alpha,t)=0$ is an elliptic equation
at $0$ on $\alpha$.  Using similar notation to the proof of
Proposition \ref{ch8s2subs1prop4},
$$d^*F_2(\alpha,t)(x)=R_t(x,\alpha(x),\nabla\alpha(x))\nabla^2\alpha(x)
+E_t(x,\alpha(x),\nabla\alpha(x)),$$ where $R_t$ and $E_t$ are
smooth functions of their arguments. If we define
$$S_{(\alpha,\,t)}(\gamma)(x)=R_t(x,\alpha(x),\nabla\alpha(x))\nabla^2\gamma(x),$$
then $S_{(\alpha,\,t)}$ is a smooth linear differential operator on
$\gamma\in C^2_{\text{loc}}(\Lambda^2_+T^*\Nhat)$. The ellipticity
of $S_{\alpha}=S_{(\alpha,0)}$ results from the coassociativity of
$\Nhat$.  Ellipticity is an \emph{open} condition so, although
$\Nhat(t)$ is not necessarily coassociative, the fact that it is
`close' to being coassociative means that $S_{(\alpha,\,t)}$ is
elliptic, as long as we shrink $\hat{\mathcal{T}}$ as necessary to
make $t$ sufficiently small.

The regularity results for $S_{(\alpha,\,t)}$ follow in the same way
as in the proof of Proposition \ref{ch8s2subs1prop4} since
$F_2(\alpha,t)$ depends smoothly on $t$ and $\Nhat(t)$ is
asymptotically coassociative near the singular points, which
validates the use of the theory from \cite[$\S$6.1.1]{Marshall}.
Recall that $L^p_{k+1,\,\lambda}\hookrightarrow C^{k,\,a}_\lambda$
where $a=1-4/p$.  Thus, if $S_{(\alpha,\,t)}(\gamma)\in
C^{k-1,\,a}_{\lambda-2}$ and $\gamma\in C^2_\lambda(\Uhat)$, then
$\gamma\in C^{k+1,\,a}_\lambda(\Uhat)$.

Since $E_0=E$ maps into $C^{k-1,\,a}_{\lambda-2}$, as argued in the
proof of Proposition \ref{ch8s2subs1prop4}, and $F_2$ depends
smoothly on $t$, $E_t$ maps into $C^{k-1,\,a}_{\lambda-2}$ for $t$
sufficiently small. Hence,
$$S_{(\alpha,\,t)}(\alpha)(x)=-E_t(x,\alpha(x),\nabla\alpha(x))
\in C^{k-1,\,a}_{\lambda-2}.$$  We deduce that $\alpha\in
C^{k+1,\,a}_\lambda$, given only that $\alpha\in C^{k,\,a}_\lambda$.
Induction gives the result.
\end{proof}

\subsection{Problem 3: moving singularities and varying $\G2$
structure}\label{ch8s2subs3}

For our final problem we consider CS deformations $N^{\prime}$ of
$N$ with the same rate and cones at $s$ singularities, but with
possibly different singular points and tangent cones there, such
that $N^{\prime}$ is coassociative under a deformation of the
$\text{G}_2$ structure on $M$.

We begin with the following. 

\begin{prop}\label{ch8s2subs3prop1}
Use the notation of Proposition \ref{ch8s2subs1prop1}.  Let
$$T=\That\cup\bigcup_{i=1}^sV_i\supseteq N.$$
By making $\That$ and $V_i$, for $i=1,\ldots,s$, smaller if
necessary, $T$ retracts onto $N$.  There exists an isomorphism
$\Xi:H^3_{\text{\emph{dR}}}(T)\rightarrow
H^3_\text{\emph{cs}}(\Nhat)$.
\end{prop}

\begin{proof} Let $[\xi]\in H^3_\text{dR}(T)$.  Since the sets
$V_i$ retract onto $\{z_i\}$ for $i=1,\ldots,s$, $\xi$ can be chosen
such that $\xi|_{V_i}=0$.  Therefore, $\xi|_{U_i\setminus\{z_i\}}=0$
which implies that the support of $\xi|_{\Nhat}$ is contained in
$K$, which is compact.  Hence $[\xi|_{\Nhat}]$ is a well-defined
element of $H^3_{\text{cs}}(\Nhat)$.  Define $\Xi$ by
$[\xi]\mapsto[\xi|_{\Nhat}]$.
We show that $\Xi$ is well-defined.  Suppose that
$\xi^\prime=\xi+d\upsilon$, for $\upsilon\in
C^{\infty}(\Lambda^2T^*T)$, such that $\xi^\prime|_{V_i}=0$ for all
$i$.  Then $d\upsilon|_{V_i}=0$ for all $i$.  Since $V_i$ retracts
onto $\{z_i\}$ we can choose $\upsilon$ such that
$\upsilon|_{V_i}=0$ without affecting $d\upsilon$ by smoothly
interpolating over $\That$.  Thus $\upsilon|_{\Nhat}$ is compactly
supported on $\Nhat$ and
$\xi|_{\Nhat}+d(\upsilon|_{\Nhat})=\xi^\prime|_{\Nhat}$.  Hence
$\Xi$ is well-defined and injective.

Any closed form on $\Nhat$ with support in $K$ can be extended
smoothly to a closed form on $T$ which vanishes on $V_i$ for all
$i$. Thus, any cohomology class in $H^3_{\text{cs}}(\Nhat)$ has a
representative $\gamma$ that can be lifted to a form $\xi$ on $T$
such that $\Xi([\xi])=[\gamma]$, which implies that $\Xi$ is
surjective.
\end{proof}

\begin{notes} The reason for this result is two-fold.  \begin{itemize}\item[(a)] The
condition that $\Xi([\varphi|_T])=0$ in $H^3_{\text{cs}}(\Nhat)$ is
implied by the coassociativity of $\Nhat$ and it forces
$[\varphi|_{\Nhat}]=0$ in $H^3_{\text{cs}}(\Nhat)$.   This is
stronger than the seemingly more natural condition of
$[\varphi|_{\Nhat}]=0$ in $H^3_{\text{dR}}(\Nhat)$, which would be
the correct requirement if $\Nhat$ were compact by the work of
McLean \cite{McLean}. \item[(b)] If a $\text{G}_2$ structure
$(\varphi^{\prime},g^{\prime})$ on $M$ is such that
$\Xi([\varphi^{\prime}|_T])\neq0$ then
$\varphi^{\prime}|_{\hat{N}^\prime}\neq0$ for any nearby deformation
$\hat{N}^{\prime}$ of $\hat{N}$, so there are \emph{no}
coassociative deformations.\end{itemize}\end{notes}

Proposition \ref{ch8s2subs3prop1} allows us to define a
distinguished family of `nearby' $\text{G}_2$ structures to
$(\varphi,g)$.

\begin{dfn}\label{ch8s2subs3dfn1}
 Let
$\hat{\mathcal{F}}$ be a small open ball about $0$ in $\R^m$ for
some $m$.  Let
$$\mathcal{F}=\{(\varphi^f,g^f)\,:\,f\in\hat{\mathcal{F}}\}$$ be a
family of torsion-free $\text{G}_2$ structures, with
$(\varphi^0,g^0)=(\varphi,g)$, such that $\Xi([\varphi^f|_T])=0$ in
$H^3_{\text{cs}}(\Nhat)$ and the map
$h_{\hat{\mathcal{F}}}:\hat{\mathcal{F}}\rightarrow\mathcal{F}$
given by $h_{\hat{\mathcal{F}}}(f)=(\varphi^f,g^f)$ is an embedding.
\end{dfn}

\begin{note} $\hat{\mathcal{F}}$ can be considered as an open subset of
$T_0\hat{\mathcal{F}}$. \end{note}

We now describe the moduli space for Problem 3.

\begin{dfn}\label{ch8s2subs3dfn2} The \emph{moduli space of deformations}
$\mathcal{M}_3(N,\lambda)$ for Problem 3 is the set of pairs
$(N^{\prime},f)$ of $f\in\hat{\mathcal{F}}$ and $N^\prime$ in
$(M,\varphi^f,g^f)$ which are CS coassociative 4-folds at
$z_1^\prime,\ldots,z_s^\prime$ with rate $\lambda$, having cone
$C_i$ and tangent cone
 $\hat{C}_i^\prime$ at $z_i^\prime$ for all $i$, such that
there exists a homeomorphism $h:N\rightarrow N^\prime$, isotopic to
the identity, with $h(z_i)=z_i^\prime$ for $i=1,\ldots,s$ and such
that $h|_{\Nhat}:\Nhat\rightarrow
N^\prime\setminus\{z_1^\prime,\ldots,z_s^\prime\}$ is a
diffeomorphism.

We have a projection map
$\pi_{\hat{\mathcal{F}}}:\mathcal{M}_3(N,\lambda)\rightarrow\hat{\mathcal{F}}$,
with $\pi_{\hat{\mathcal{F}}}(N^\prime,f)=f$, whose fibres
$\pi_{\hat{\mathcal{F}}}^{-1}(f)$ are equal to the moduli space for
Problem 2 defined using the $\text{G}_2$ structure
$(\varphi^f,g^f)$.
\end{dfn}

We must adapt our translation space from Problem 2 to incorporate
the varying $\text{G}_2$ structure.

\begin{dfn}\label{ch8s2subs3dfn5}
 Use the notation of Definitions \ref{ch8s2subs2dfn2} and \ref{ch8s2subs3dfn1}.
For $f\in\hat{\mathcal{F}}$ and $\mathbf{z}^\prime\in B$ let
$I^f(\mathbf{z}^\prime)$ denote the family of choices of $s$-tuples
$\boldsymbol{\zeta}^\prime=(\zeta_1^\prime,\ldots,\zeta_s^\prime)$
of isomorphisms $\zeta_i^\prime:\R^7\rightarrow T_{z_i^\prime}M$
identifying $(\varphi_0,g_0)$ with $(\varphi^f|_{T_{z_i^\prime}M},
g^f|_{T_{z_i^\prime}M})$.

 The \emph{translation space}
corresponding to $\hat{\mathcal{F}}$ is
$$\mathcal{T}^{\hat{\mathcal{F}}}=\{(\mathbf{z}^\prime,\boldsymbol{\zeta}^\prime,f)\,:\,
\mathbf{z}^\prime\in
B,\,f\in\hat{\mathcal{F}},\,\boldsymbol{\zeta}^\prime\in
I^f(\mathbf{z}^\prime)\}.$$ It is a principal $\text{G}_2^s$ bundle
over $B\times\hat{\mathcal{F}}$.

There is a natural free action of $\text{H}$ on
$\mathcal{T}^{\hat{\mathcal{F}}}$ and hence an $\text{H}$-orbit
through $(\mathbf{z},\boldsymbol{\zeta},0)$.
 Therefore, we may embed
$\hat{\mathcal{T}}\times\hat{\mathcal{F}}$ into
$\mathcal{T}^{\hat{\mathcal{F}}}$ by
$h_{\hat{\mathcal{T}}\times\hat{\mathcal{F}}}:(t,f)\mapsto
(\mathbf{z}(t,f),\boldsymbol{\zeta}(t,f),f)$ such that
$h_{\hat{\mathcal{T}}\times\hat{\mathcal{F}}}(\hat{\mathcal{T}}\times\hat{\mathcal{F}})$
is transverse to this $\text{H}$-orbit,\linebreak
$h_{\hat{\mathcal{T}}\times\hat{\mathcal{F}}}(t,0)=h_{\hat{\mathcal{T}}}(t)$
for all $t$ and $\mathbf{z}(0,f)=\mathbf{z}$ for all $f$.
\end{dfn}

Use the notation introduced before Proposition
\ref{ch8s2subs2prop1}. Extend the $\G2$ coordinate system near
$z_1,\ldots,z_s$ used to define $N$ to a smooth family of
$\text{G}_2$ coordinate systems
$$\left\{\{\chi_i(t,f):B(0;\eta)\rightarrow V_i(t,f)\,:\,i=1,\ldots,s\}\,:\,(t,f)
\in\hat{\mathcal{T}}\times\hat{\mathcal{F}}\right\}$$ such that
$V_i(t,f)$ is an open set in $M$ containing $z_i(t,f)$,
$\chi_i(t,f)(0)=z_i(t,f)$, $d\chi_i(t,f)|_0=\zeta_i(t,f)$,
$\chi_i(t,0)=\chi_i(t)$, $V_i(0,f)=V_i$ and $V_i(t,0)=V_i(t)$ for
$i=1,\ldots,s$. We state the analogue of Proposition
\ref{ch8s2subs2prop1}.

\begin{prop}\label{ch8s2subs3prop2}Use the notation of Propositions
\ref{ch8s2subs1prop1} and \ref{ch8s2subs2prop1} and Definition
\ref{ch8s2subs3dfn5}.
\begin{itemize}\item[(a)] There exists a family
$\mathcal{N}^{\hat{\mathcal{F}}}=\{N(t,f):(t,f)\in\hat{\mathcal{T}}\times\hat
{\mathcal{F}}\}$ of CS 4-folds in $M$, with $N(0,f)=N$ and
$N(t,0)=N(t)$, such that $N(t,f)$ has singularities at
$z_1(t,f),\ldots,z_s(t,f)$ with rate $\lambda$, cones
$C_1,\ldots,C_s$ and tangent cones
$\hat{C}_1(t,f),\ldots,\hat{C}_s(t,f)$ defined by
$\hat{C}_i(t,f)=(\zeta_i(t,f)\circ\iota_i)(C_i)$.
\item[(b)]  Let
$\hat{N}(t,f)=N(t,f)\setminus\{z_1(t,f),\ldots,z_s(t,f)\}$ and write
$$N(t,f)=K(t,f)\sqcup\bigsqcup_{i=1}^sU_i(t,f)$$
where $K(t,f)$ is compact and
$U_i(t,f)\setminus\{z_i(t,f)\}\cong(0,\epsilon)\times\Sigma_i$ for
all $i$, in the obvious way, ensuring that $K(0,f)=K$,
$K(t,0)=K(t)$, $U_i(0,f)=U_i$ and $U_i(t,0)=U_i(t)$.  For
$(t,f)\in\hat{\mathcal{T}}\times\hat{\mathcal{F}}$, there exist open
sets $\hat{T}(t,f)\subseteq M$ containing $\hat{N}(t,f)$ and
diffeomorphisms $\delta(t,f):\hat{U}\rightarrow\hat{T}(t,f)$ taking
the zero section to $\hat{N}(t,f)$, varying smoothly in $t$ and $f$,
with $\hat{T}(0,f)=\hat{T}$, $\hat{T}(t,0)=\hat{T}(t)$ and
$\delta(t,0)=\delta(t)$. Moreover, $\hat{T}(t,f)$ can be chosen to
grow with order $O(\rho)$ as $\rho\rightarrow 0$ and $\delta(t,f)$
is compatible with the identifications
$U_i(t,f)\setminus\{z_i(t,f)\}\cong (0,\epsilon)\times\Sigma_i$ for
$i=1,\ldots,s$.
\end{itemize}
\end{prop}

\noindent The proof is almost identical to that of Proposition
\ref{ch8s2subs2prop1} and so we omit it.  The compatibility
conditions on $\delta(t,f)$ are given by similar commutative
diagrams to \eq{ch8s2subs2eq1}.

\begin{remark} $\delta(t,f)$ is not required to be compatible with the
isomorphism $\nu(\Nhat)\cong\Lambda^2_+T^*\Nhat$ for $(t,f)\neq
(0,0)$.\end{remark}

\medskip

We proceed by defining our final deformation map.

\begin{dfn}\label{ch8s2subs3dfn3}
Use the notation of Proposition \ref{ch8s2subs3prop2}.  Let
$\Gamma_\alpha$ be the graph of $\alpha\in C^1_{\text{loc}}(\Uhat)$
and let $\pi_\alpha:\Nhat\rightarrow\Gamma_\alpha$ be given by
$\pi_\alpha(x)=(x,\alpha(x))$.  For
$(t,f)\in\hat{\mathcal{T}}\times\hat{\mathcal{F}}$, let
$f_\alpha(t,f)=\delta(t,f)\circ\pi_\alpha$ and let
$\Nhat_\alpha(t,f)=f_\alpha(t,f)(\Nhat)$.  Define $F_3$ from
$C^1_{\text{loc}}(\Uhat)\times\hat{\mathcal{T}}\times\hat{\mathcal{F}}$
to $C^0_{\text{loc}}(\Lambda^3T^*\Nhat)$ by:
$$F_3(\alpha,t,f)=f_\alpha(t,f)^*\left(\varphi^f|_{\Nhat_\alpha(t,f)}\right).$$
The linearisation of $F_3$ at $(0,0,0)$ acts as
$$dF_3|_{(0,0,0)}:(\alpha,t,f)\longmapsto d\alpha+L_2(t)+L_3(f),$$
where $\alpha\in C^1_{\text{loc}}(\Lambda^2_+T^*\Nhat)$, $(t,f)\in
T_0\hat{\mathcal{T}}\oplus T_0\hat{\mathcal{F}}$, $L_2$ is given in
Definition \ref{ch8s2subs2dfn3} and $L_3$ is a linear map into the
space of smooth exact 3-forms on $\Nhat$ by the condition imposed on
$\varphi^f$ in Definition \ref{ch8s2subs3dfn1}.
\end{dfn}

\begin{note} $F_3(\alpha,t,0)=F_2(\alpha,t)$ as given in Definition
\ref{ch8s2subs2dfn3}. \end{note}

Now, $\text{Ker}\,F_3$ corresponds to choices of $\Nhat_\alpha(t,f)$
which are coassociative with respect to $(\varphi^f,g^f)$. The next
result is then clear from considering the proof of Proposition
\ref{ch8s2subs2prop2}.

\begin{prop}\label{ch8s2subs3prop3}
The moduli space of deformations for Problem 3 is locally
homeomorphic to $\text{\emph{Ker}}\,F_3=\{(\alpha,t,f)\in
C^{\infty}_{\lambda}(\Uhat)\times\hat{\mathcal{T}}
\times\hat{\mathcal{F}}\,:\, F_3(\alpha,t,f)=0\}$.
\end{prop}

We again have an associated map to our deformation map.

\begin{dfn}\label{ch8s2subs3dfn4} Define $G_3:C_{\text{loc}}^{1}(\Uhat)\times
C_{\text{loc}}^{1}(\Lambda^4T^*\Nhat)\times\hat{\mathcal{T}}\times\hat{\mathcal{F}}
\rightarrow C^{0}_{\text{loc}}(\Lambda^3T^*\Nhat)$ by:
\begin{equation*}
G_3(\alpha,\beta,t,f)=F_3(\alpha,t,f)+d^*\beta. \end{equation*} Then
$dG_3|_{(0,0,0,0)}:(\alpha,\beta,t,f)\longmapsto
d\alpha+d^*\beta+L_2(t)+L_3(f)$, in the notation of Definition
\ref{ch8s2subs3dfn3}.
\end{dfn}

\noindent The next result is analogous to Propositions
\ref{ch8s2subs1prop3} and \ref{ch8s2subs2prop3} and may be
immediately deduced from the exactness of $F_3(\alpha,t,f)$, which
follows from the condition imposed on $\varphi^f$ in Definition
\ref{ch8s2subs3dfn1}.

\begin{prop}\label{ch8s2subs3prop4}
$$\text{\emph{Ker}}\,F_3\cong\{(\alpha,\beta,t,f)\in
C^{\infty}_{\lambda}(\Uhat)\times
C^{\infty}_{\lambda}(\Lambda^4T^*\Nhat)\times\hat{\mathcal{T}}
\times\hat{\mathcal{F}}\,:\,G_3(\alpha,\beta,t,f)=0\}.$$
\end{prop}

The argument used to prove the regularity result Proposition
\ref{ch8s2subs2prop4} is easily generalised to the map $G_3$, so we
end the section with the following.

\begin{prop}\label{ch8s2subs3prop5}
Let $(\alpha,\beta,t,f)\in L^p_{k+1,\,\lambda}(\Uhat)\times
L^p_{k+1,\,\lambda}(\Lambda^4T^*\Nhat)\times\hat{\mathcal{T}}\times\hat{\mathcal{F}}$,
where $p>4$ and $k\geq 2$. If\/ $G_3(\alpha,\beta,t,f)=0$ and
$\|\alpha\|_{C^1_1}$, $t$ and $f$ are sufficiently small,
$(\alpha,\beta)\in C^{\infty}_\lambda(\Uhat)\times
C^{\infty}_\lambda(\Lambda^4T^*\Nhat)$.
\end{prop}

\section{The deformation and obstruction spaces}
\label{ch8s3}

In this section we describe the infinitesimal deformation and
obstruction spaces for each of our problems and show in each
scenario that, if the obstruction space is zero, we get a
\emph{smooth} moduli space of deformations. We recollect the common
notation introduced at the start of $\S$\ref{ch8s2}.  In addition,
fix some $p>4$ and integer $k\geq 2$.

\subsection{Problem 1}\label{ch8s3subs1}

Recall the maps $F_1$ and $G_1$ given in Definitions
\ref{ch8s2subs1dfn2} and \ref{ch8s2subs1dfn3} respectively.  Their
kernels give a local description for the moduli space
$\mathcal{M}_1(N,\lambda)$ by Propositions \ref{ch8s2subs1prop2} and
\ref{ch8s2subs1prop3}.  Therefore the kernels of $dF_1|_0$ and
$dG_1|_{(0,0)}$ describe the \emph{infinitesimal deformations}.

\begin{dfn}\label{ch8s3subs1dfn1}
The \emph{infinitesimal deformation space} for Problem 1 is
\begin{align*}
\mathcal{I}_1(N,\lambda)&=\{\alpha\in
C^{\infty}_{\lambda}(\Lambda^2_+ T^*\Nhat)\,:\,d\alpha=0\}\\
&\cong\{(\alpha,\beta)\in C^{\infty}_{\lambda} (\Lambda^2_+T^*\Nhat
\oplus \Lambda^4T^*\Nhat)\,:\,d\alpha+d^*\beta=0\}.
\end{align*}
The equivalence of the spaces follows by Proposition
\ref{ch8s2subs1prop3} or, more simply, by the Maximum Principle for
harmonic functions.

Using Proposition \ref{ch8s2subs1prop4},
$$\mathcal{I}_1(N,\lambda)\cong\{(\alpha,\beta)\in
L^{p}_{k+1,\,\lambda} (\Lambda^2_+T^*\Nhat
\oplus\Lambda^4T^*\Nhat)\,:\,d\alpha+d^*\beta=0\}.$$ Therefore,
$\mathcal{I}_1(N,\lambda)$ is finite-dimensional.
\end{dfn}

We turn to possible \emph{obstructions} to the deformation theory
and start with the following.

\begin{prop}\label{ch8s3subs1prop1}
The map $F_1$ 
takes $L^p_{k+1,\,\lambda}(\Uhat)$ into
$d(L^p_{k+1,\,\lambda}(\Lambda^2T^*\Nhat))$.
\end{prop}

\begin{proof}
Let $\alpha\in L^p_{k+1,\,\lambda}(\Uhat)$ and let $T$ be as in
Proposition \ref{ch8s2subs3prop1}.  As noted after that proposition,
$[\varphi|_T]=0$ in $H^3_{\text{dR}}(T)$ and hence $\varphi|_T$ is
exact.  Thus, $\varphi|_T=d\psi$ for some $\psi\in
C^{\infty}(\Lambda^2T^*T)$. However, we want to select $\psi$ in a
particular way near the singularities. On $B(0;\eta)\subseteq\R^7$,
for each $i=1,\ldots,s$,
$$\chi_i^*(\varphi)=\varphi_0+O(r_i).$$  If
$v$ is the dilation vector field on $\R^7$, given in coordinates
$(x_1,\ldots,x_7)$ by
\begin{equation*}
v=x_1\frac{\partial}{\partial
x_1}+\ldots+x_7\frac{\partial}{\partial x_7}\,,\end{equation*} we
can choose $\psi$ to satisfy
$$\chi_i^*(\psi)=\frac{1}{3}(v\cdot\varphi_0)+O(r_i^2)$$
on $V_i$, since $d(v\cdot\varphi_0)=3\varphi_0$,
 then extend $\psi$ smoothly to a form on $T$
  such that
$d\psi=\varphi|_T$. Note that
$$(v\cdot\varphi_0)|_{\iota_i(C_i)}=v\cdot(\varphi_0|_{\iota_i(C_i)})=0$$
 as $v\in T(\iota_i(C_i))$. Hence
$\chi^*_i(\psi)=O(r_i^2)$ on $\iota_i(C_i)$, for all $i$, and
similar results hold for the derivatives of $\psi$.  Define
$$H_1(\alpha)=f_{\alpha}^*\left(\psi|_{\Nhat_\alpha}\right)$$
so that $F_1(\alpha)=d(H_1(\alpha))$.  Note that
$\chi_i^*(\psi)|_{\iota_i(C_i)}=O(r_i^2)$ is dominated by
$O(r_i^\lambda)$ terms as $r_i\rightarrow0$ since $\lambda<2$.
Further, $f_{\alpha}^*(\psi|_{\Nhat_\alpha})$ has the same growth as
$\chi_i^*(\psi)|_{(\Phi_\alpha)_i((0,\epsilon)\times\Sigma_i)}$ as
$r_i\rightarrow 0$, using the notation preceding Proposition
\ref{ch8s2subs1prop2}. However,
$$\chi_i^*(\psi)|_{(\Phi_\alpha)_i((0,\epsilon)\times\Sigma_i)}=
\chi_i^*(\psi)|_{((\Phi_\alpha)_i-\iota_i)((0,\epsilon)\times\Sigma_i)}+
\chi_i^*(\psi)|_{\iota_i((0,\epsilon)\times\Sigma_i)}.$$ The first
term on the right-hand side depends on $|(\Phi_\alpha)_i-\iota_i|$
and hence is $O(r_i^\lambda)$ as $r_i\rightarrow 0$.  This dominates
the second term by our observation above.  Hence, $H_1(\alpha)\in
L^p_{k,\,\lambda}$ because $H_1$ depends on $\alpha$ and
$\nabla\alpha$.
Note that $H_1(\alpha)$ has one degree of differentiability less
than expected.

Recalling that $\lambda\notin\mathcal{D}$, we deduce that
$F_1(\alpha)$ lies in $d(L^p_{k,\,\lambda}(\Lambda^2T^*\Nhat))$ and
hence is $L^2$-orthogonal to elements of the kernel of
$$d+d^*:L^q_{l+1,\,-3-\lambda}(\Lambda^3T^*\Nhat)
\rightarrow
L^q_{l,\,-4-\lambda}(\Lambda^2T^*\Nhat\oplus\Lambda^4T^*\Nhat),$$
where $q>1$ such that $1/p+1/q=1$.  We show that
$$d(L^p_{k,\,\lambda}(\Lambda^2T^*\Nhat))\oplus d^*(L^p_{k,\,\lambda}
(\Lambda^4T^*\Nhat))\subseteq
L^p_{k-1,\,\lambda-1}(\Lambda^3T^*\Nhat)$$ is \emph{characterised}
as the subspace which is $L^2$-orthogonal to this kernel.

Consider
$$d+d^*:L^p_{k,\,\lambda}(\Lambda^\text{even}T^*\Nhat)\rightarrow
L^p_{k-1,\,\lambda-1}(\Lambda^{\text{odd}}T^*\Nhat).$$  This
elliptic map has image which comprises precisely of those elements
of $L^p_{k-1,\,\lambda-1}(\Lambda^{\text{odd}}T^*\Nhat)$ which are
$L^2$-orthogonal to the kernel $\mathcal{K}$ of
$$d+d^*:L^q_{l+1,\,-3-\lambda}(\Lambda^{\text{odd}}T^*\Nhat)
\rightarrow L^q_{l,\,-4-\lambda}(\Lambda^\text{even}T^*\Nhat).$$ The
space $\mathcal{K}$ can be written as the direct sum
$\mathcal{K}=\mathcal{K}^1\oplus\mathcal{K}^3\oplus\mathcal{K}^m$,
where
$$\mathcal{K}^j=\mathcal{K}\cap L^q_{l+1,\,-3-\lambda}(\Lambda^{j}T^*\Nhat)$$
for $j=1$ and $3$ and $\mathcal{K}^m$ is some transverse subspace.
Then
$$d(L^p_{k,\,\lambda}(\Lambda^2T^*\Nhat))\oplus d^*(L^p_{k,\,\lambda}
(\Lambda^4T^*\Nhat))=\{\alpha_3\,:\,\exists\,\alpha_1\;\text{such
that}\;(\alpha_1,\alpha_3)\in\mathcal{K}^{\perp}\},$$ where we take
the orthogonal complement in $L^p_{k-1,\,\lambda-1}$. Note that the
projection $\pi_1(\mathcal{K}^m)$ of $\mathcal{K}^m$ onto the space
of 1-forms must meet $\mathcal{K}^1$ in the zero form since, if
$(\alpha_1,\alpha_3)\in\mathcal{K}^m$ and $\alpha_1\in\mathcal{K}^1$
then $\alpha_3\in\mathcal{K}^3$, which contradicts the direct sum
decomposition of $\mathcal{K}$. Therefore, $\pi_1(\mathcal{K}^m)$
and $\mathcal{K}^1$ are transverse finite-dimensional subspaces of
$L^q_{l+1,\,-3-\lambda}(\Lambda^1T^*\Nhat)$.  Hence, there exists a
space $\mathcal{A}$ of smooth compactly supported 1-forms on $\Nhat$
which is $L^2$-orthogonal to $\mathcal{K}^1$ and such that
$\mathcal{A}\times\mathcal{K}^m\rightarrow\R$ given by
$(\gamma,\xi)\mapsto(\gamma,0)\cdot\xi$
is a dual pairing.  If $\alpha_3\in
L^p_{k-1,\,\lambda-1}(\Lambda^3T^*\Nhat)$ such that $\alpha_3\in
(\mathcal{K}^3)^{\perp}$, there exists a unique
$\alpha_1\in\mathcal{A}$ such that
$(\alpha_1,0)\cdot\xi=-(0,\alpha_3)\cdot\xi$ for all
$\xi\in\mathcal{K}^m$, which implies that
$(\alpha_1,\alpha_3)\in(\mathcal{K}^m)^{\perp}$.  We conclude that
\begin{align*}
(\mathcal{K}^3)^{\perp}&=\{\alpha_3\in(\mathcal{K}^3)^{\perp}\,:
\,\exists\, \alpha_1\in(\mathcal{K}^1)^\perp\;\text{such
that}\;(\alpha_1,\alpha_3)\in(\mathcal{K}^m)^\perp\}\\
&= \{\alpha_3\,:\,\exists\,\alpha_1\;\text{such
that}\;(\alpha_1,\alpha_3)\in\mathcal{K}^{\perp}\}\\
&=d(L^p_{k,\,\lambda}(\Lambda^2T^*\Nhat))\oplus
d^*(L^p_{k,\,\lambda} (\Lambda^4T^*\Nhat))\subseteq
L^p_{k-1,\,\lambda-1}(\Lambda^3T^*\Nhat).
\end{align*}

  However,
$\mathcal{K}^3$ is independent of $k$, and hence $F_1(\alpha)$ must
lie in the image of $d+d^*$ from $L^p_{k+1,\,\lambda}$, since
$F_1(\alpha)$ lies in $L^p_{k,\,\lambda-1}$.  We may thus write
$F_1(\alpha)=d\gamma+d^*\beta$ for some $\gamma\in
L^p_{k+1,\,\lambda}(\Lambda^2T^*\Nhat)$ and $\beta\in
L^p_{k+1,\,\lambda}(\Lambda^4T^*\Nhat)$. Moreover, $dd^*\beta=0$ and
so $\beta$ is harmonic and $O(\rho^\lambda)$ as $\rho\rightarrow 0$.
By the Maximum Principle, (noting that $*\beta$ is a harmonic
function on $\Nhat$), $\beta=0$.  The proposition is thus proved.
\end{proof}

 We deduce from Propositions \ref{ch8s2subs1prop2},
\ref{ch8s2subs1prop3}, \ref{ch8s2subs1prop4} and
\ref{ch8s3subs1prop1} that $\mathcal{M}_1(N,\lambda)$ is locally
homeomorphic to the kernel of
$$G_1:L^p_{k+1,\,\lambda}(\Uhat)\times L^p_{k+1,\,\lambda}(\Lambda^4T^*\Nhat)
\rightarrow d(L^p_{k+1,\,\lambda}(\Lambda^2T^*\Nhat))\oplus d^*(
L^p_{k+1,\,\lambda}(\Lambda^4T^*\Nhat)).$$  Therefore, our
deformation theory will be obstructed if and only if the map
$$d:L^p_{k+1,\,\lambda}(\Lambda^2_+T^*\Nhat)\rightarrow
d(L^p_{k+1,\,\lambda}(\Lambda^2T^*\Nhat))$$ is \emph{not}
surjective.  This leads us to the next result and definition.

\begin{prop}\label{ch8s3subs1prop2}
There exists a finite-dimensional subspace
$\mathcal{O}_1(N,\lambda)$ of\/\\
$L^p_{k,\,\lambda-1}(\Lambda^3T^*\Nhat)$ such that
$$d(L^p_{k+1,\,\lambda}(\Lambda^2T^*\Nhat))=d(L^p_{k+1,\,\lambda}
(\Lambda^2_+T^*\Nhat))\oplus\mathcal{O}_1(N,\lambda).$$
\end{prop}

\begin{proof} The Fredholmness of $d+d^*$
implies that the images of
$L^p_{k+1,\,\lambda}(\Lambda^2_+T^*\Nhat)\oplus
L^p_{k+1,\,\lambda}(\Lambda^4T^*\Nhat)$ and
$L^p_{k+1,\,\lambda}(\Lambda^2T^*\Nhat)\oplus
L^p_{k+1,\,\lambda}(\Lambda^4T^*\Nhat)$ under $d+d^*$ are both
closed and have finite codimension in
$L^p_{k,\,\lambda-1}(\Lambda^3T^*\Nhat)$.  Since
\begin{align*} \{0\}&=d(L^p_{k+1,\,\lambda}(\Lambda^2T^*\Nhat))\cap
d^*(L^p_{k+1,\,\lambda}(\Lambda^4T^*\Nhat))\\
&=d(L^p_{k+1,\,\lambda}(\Lambda^2_+T^*\Nhat))\cap
d^*(L^p_{k+1,\,\lambda}(\Lambda^4T^*\Nhat))\end{align*} by the
Maximum Principle, we deduce that
$$d(L^p_{k+1,\,\lambda}(\Lambda^2_+T^*\Nhat))\quad\text{and}\quad
d(L^p_{k+1,\,\lambda}(\Lambda^2T^*\Nhat))$$ are both closed and that
the former has finite codimension in the latter.  Thus,
$\mathcal{O}_1(N,\lambda)$ can be chosen as stated.
\end{proof}

\begin{dfn}\label{ch8s3subs1dfn2}
The \emph{obstruction space} for Problem 1 is
$$\mathcal{O}_1(N,\lambda)\cong\frac{d(L^p_{k+1,\,\lambda}
(\Lambda^2T^*\Nhat))}{d(L^p_{k+1,\,\lambda}(\Lambda^2_+T^*\Nhat))}\,.$$
\end{dfn}


We proceed as follows.  Define
\begin{align*}
U_1&=L^p_{k+1,\,\lambda}(\Uhat)\times L^p_{k+1,\,\lambda}
(\Lambda^4T^*\Nhat),\\
X_1&=L^p_{k+1,\,\lambda}(\Lambda^2_+T^*\Nhat\oplus\Lambda^4T^*\Nhat),\\
Y_1&=\mathcal{O}_1(N,\lambda)\subseteq
L^p_{k,\,\lambda-1}(\Lambda^3T^*\Nhat)\;\,\text{and}\\
Z_1&=d(L^p_{k+1,\,\lambda}(\Lambda^2T^*\Nhat))\oplus d^*(
L^p_{k+1,\,\lambda}(\Lambda^4T^*\Nhat)).
\end{align*}
Then $X_1$, $Y_1$ and $Z_1$ are Banach spaces and $U_1$ is an open
neighbourhood of $(0,0)$ in $X_1$ because
$L^p_{k+1,\,\lambda}\hookrightarrow C^0_\lambda$ by Theorem
\ref{ch6s2thm1} and $\Uhat$ grows with order $O(\rho)$ as
$\rho\rightarrow0$ by Proposition \ref{ch8s2subs1prop1}. Thus,
$W_1=U_1\times Y_1$ is an open neighbourhood of $(0,0,0)$ in
$X_1\times Y_1$. Define $\mathcal{G}_1:W_1\rightarrow Z_1$ by:
$$\mathcal{G}_1(\alpha,\beta,\gamma)=G_1(\alpha,\beta)+\gamma.$$ Then
$\mathcal{G}_1$ is well-defined by Propositions
\ref{ch8s3subs1prop1} and \ref{ch8s3subs1prop2} and its derivative
at $(0,0,0)$ acts from $X_1\times Y_1$ to $Z_1$ as
\begin{align*} 
d\mathcal{G}_1|_{(0,0,0)}:(\alpha,\beta,\gamma)\longmapsto
d\alpha+d^*\beta+\gamma.
\end{align*}
Clearly, $d\mathcal{G}_1|_{(0,0,0)}$ is surjective by construction
and its kernel, using the fact that $(d+d^*)(X_1)\cap Y_1=\{0\}$, is
given by:
\begin{align*}
\text{Ker}\,d\mathcal{G}_1|_{(0,0,0)}&=\{(\alpha,\beta,\gamma)\in
X_1\times Y_1:d\alpha+d^*\beta+\gamma=0\}\\
&\cong\{(\alpha,\beta)\in
X_1:d\alpha+d^*\beta=0\}\cong\mathcal{I}_1(N,\lambda).
\end{align*}
The conclusion, by implementing the Implicit Function Theorem for
Banach spaces (Theorem \ref{ch6s2thm2}), is that $\text{Ker}\,
\mathcal{G}_1$ is a smooth manifold near zero which may be
identified with an open neighbourhood
$\hat{\mathcal{M}}_1(N,\lambda)$ of $0$ in
$\mathcal{I}_1(N,\lambda)$.  Formally, if we write
$X_1=\mathcal{I}_1(N,\lambda)\oplus A$ for some closed subspace $A$
of $X_1$, there exist open sets
$\hat{\mathcal{M}}_1(N,\lambda)\subseteq \mathcal{I}_1(N,\lambda)$,
$V_A\subseteq A$, $V_Y\subseteq Y_1$, all containing $0$, with
$\hat{\mathcal{M}}_1(N,\lambda)\times V_A\subseteq U_1$, and smooth
maps $\mathcal{V}_A:\hat{\mathcal{M}}_1(N,\lambda)\rightarrow V_A$
and $\mathcal{V}_Y:\hat{\mathcal{M}}_1(N,\lambda)\rightarrow V_Y$
such that
$$\text{Ker}\,\mathcal{G}_1\cap(\hat{\mathcal{M}}_1(N,\lambda)\times V_A\times V_Y)
=\{(x,\mathcal{V}_A(x),\mathcal{V}_Y(x))\,:\,x\in
\hat{\mathcal{M}}_1(N,\lambda)\}.$$ If we define a smooth map
$\pi_1:\hat{\mathcal{M}}_1(N,\lambda)\rightarrow\mathcal{O}_1(N,\lambda)$
by $\pi_1(x)=\mathcal{V}_Y(x)$, the moduli space
$\mathcal{M}_1(N,\lambda)$ near $N$ is locally homeomorphic to the
kernel of $\pi_1$ near $0$.  We can think of $\pi_1$ as a map on an
open neighbourhood of $(0,0,0)$ in $\text{Ker}\,\mathcal{G}_1$ which
projects onto the obstruction space.  We write these results as a
theorem.

\begin{thm}\label{ch8s3subs1thm1} Use the notation of Definitions
\ref{ch8s2subs1dfn1}, \ref{ch8s3subs1dfn1} and \ref{ch8s3subs1dfn2}.
There exists a smooth manifold\/ $\hat{\mathcal{M}}_1(N,\lambda)$,
which is an open neighbourhood of\/ $0$ in
$\mathcal{I}_1(N,\lambda)$, and a smooth map\/
$\pi_1:\hat{\mathcal{M}}_1(N,\lambda)\rightarrow\mathcal{O}_1(N,\lambda)$,
with $\pi_1(0)=0$, such that an open neighbourhood of\/ $0$ in
$\text{\emph{Ker}}\,\pi_1$ is homeomorphic to an open neighbourhood
of\/ $N$ in $\mathcal{M}_1(N,\lambda)$.
\end{thm}

\noindent We deduce from this theorem that, if the obstruction space
is zero, the moduli space is a smooth manifold near $N$ of dimension
equal to that of the infinitesimal deformation space.  We expect the
obstruction space to be zero for generic choices of $N$ and the
$\text{G}_2$ structure on $M$.

\subsection{Problem 2}
\label{ch8s3subs2}

Recall the notation introduced in Definitions \ref{ch8s2subs2dfn2},
\ref{ch8s2subs2dfn3} and \ref{ch8s2subs2dfn4}.  We begin by defining
the infinitesimal deformation space for this problem.

\begin{dfn}\label{ch8s3subs2dfn1}
The \emph{infinitesimal deformation space} for Problem 2 is
\begin{align*}
\mathcal{I}_2(N,\lambda)&=\{(\alpha,t)\in
C^{\infty}_\lambda(\Lambda^2_+T^*\Nhat)\oplus
T_0\hat{\mathcal{T}}\,:\,d\alpha+L_2(t)=0\}\\
&\cong\{(\alpha,\beta,t)\in
C^{\infty}_\lambda(\Lambda^2_+T^*\Nhat\oplus\Lambda^4T^*\Nhat)\oplus
T_0\hat{\mathcal{T}}:d\alpha+d^*\beta+L_2(t)=0\}.
\end{align*}
\noindent  The equivalence in the definition follows from
Proposition \ref{ch8s2subs2prop3} or from the observation that
$d\alpha+L_2(t)$ is exact and so $\beta=0$ by the Maximum Principle.

By Proposition \ref{ch8s2subs2prop4},
$$\mathcal{I}_2(N,\lambda)\cong\{(\alpha,\beta,t)\in
L^{p}_{k+1,\,\lambda} (\Lambda^2_+T^*\Nhat
\oplus\Lambda^4T^*\Nhat)\oplus
T_0\hat{\mathcal{T}}\,:\,d\alpha+d^*\beta+L_2(t)=0\}.$$ Therefore,
$\mathcal{I}_2(N,\lambda)$ is finite-dimensional.
\end{dfn}

\begin{note}
There is a subspace of $\mathcal{I}_2(N,\lambda)$ which is
isomorphic to $\mathcal{I}_1(N,\lambda)$.
\end{note}

To start our consideration of obstructions, we have the
generalisation of Proposition \ref{ch8s3subs1prop1}.

\begin{prop}\label{ch8s3subs2prop1}
The map $F_2$ 
takes $L^p_{k+1,\,\lambda}(\Uhat)\times\hat{\mathcal{T}}$ into
$d(L^p_{k+1,\,\lambda}(\Lambda^2T^*\Nhat))$.
\end{prop}

\begin{proof} Use the notation from Proposition
\ref{ch8s2subs2prop1} and its proof and from the proof of
Proposition \ref{ch8s3subs1prop1}.  Recall that we have an open set
$T\supseteq\That$ in $M$ containing $N$, which retracts onto $N$,
and $\psi\in C^{\infty}(\Lambda^2T^*T)$ such that
$d\psi=\varphi|_T$.  We may similarly construct open sets
$T(t)\supseteq\That(t)$ in $M$, with $T(0)=T$, which contain $N(t)$
and retract onto it, varying smoothly with $t$.  We also have
$\psi(t)\in C^{\infty}(\Lambda^2T^*T(t))$, with $\psi(0)=\psi$, such
that $d\psi(t)=\varphi|_{T(t)}$, using the fact that $\varphi$ is
exact on $N(t)$.  Again, the $\psi(t)$ vary smoothly with $t$.
Formally, let
$$T(t)=\hat{T}(t)\cup\bigcup_{i=1}^sV_i(t).$$
By making $\hat{T}(t)$ and $V_i(t)$ smaller if necessary, $T(t)$
will be an open set as stated.  We may choose $\psi(t)$ such that
$$\chi_i(t)^*(\psi(t))=\frac{1}{3}\,(v\cdot\varphi_0)+O(r_i^2)$$
on $V_i(t)$ and then extend smoothly to a form $\psi(t)$ on $T(t)$
as required.  Define
$$H_2(\alpha,t)=f_\alpha(t)^*\left(\psi(t)|_{\Nhat_\alpha(t)}\right).$$
Then $d(H_2(\alpha,t))=F_2(\alpha,t)$.  Moreover, by the same
reasoning that $H_1(\alpha)\in L^p_{k,\,\lambda}$ in the proof of
Proposition \ref{ch8s3subs1prop1}, $H_2(\alpha,t)$ lies in
$L^p_{k,\,\lambda}$. Therefore, $F_2(\alpha,t)$ lies in
$d(L^p_{k,\,\lambda}(\Lambda^2T^*\Nhat))$. However, because
$F_2(\alpha,t)\in L^p_{k,\,\lambda-1}(\Lambda^3T^*\Nhat)$, the
argument at the end of the proof of Proposition
\ref{ch8s3subs1prop1} implies that $F_2(\alpha,t)\in
d(L^p_{k+1,\,\lambda}(\Lambda^2T^*\Nhat))$ as required.
\end{proof}

We now define the obstruction space.

\begin{dfn}\label{ch8s3subs2dfn2}
From Propositions \ref{ch8s3subs1prop2} and \ref{ch8s3subs2prop1},
since $L_2$ is a linear map on a finite-dimensional vector space,
there exists a finite-dimensional subspace
$\mathcal{O}_2(N,\lambda)$ of
$L^p_{k,\,\lambda-1}(\Lambda^3T^*\Nhat)$ such that
$$
d(L^p_{k+1,\,\lambda}(\Lambda^2T^*\Nhat))=(d(L^p_{k+1,\,\lambda}(\Lambda^2_+T^*\Nhat))+
L_2(T_0\hat{\mathcal{T}}))\oplus\mathcal{O}_2(N,\lambda).$$ We
define $\mathcal{O}_2(N,\lambda)$ to be the \emph{obstruction space}
for Problem 2.
\end{dfn}

\begin{note} $\mathcal{O}_2(N,\lambda)$ may be chosen to be contained in
$\mathcal{O}_1(N,\lambda)$.\end{note}

Following the scheme for Problem 1, we let
\begin{align*}
U_2&=L^p_{k+1,\,\lambda}(\Uhat)\times
L^p_{k+1,\,\lambda}(\Lambda^4T^*\Nhat)\times\hat{\mathcal{T}},\\
X_2&=L^p_{k+1,\,\lambda}(\Lambda^2_+T^*\Nhat\oplus
\Lambda^4T^*\Nhat)\oplus T_0\hat{\mathcal{T}},\\
Y_2&=\mathcal{O}_2(N,\lambda)\subseteq
L^p_{k,\,\lambda-1}(\Lambda^3T^*\Nhat)\;\,\text{and}\\
Z_2&=d(L^p_{k+1,\,\lambda}(\Lambda^2T^*\Nhat))\oplus
d^*(L^p_{k+1,\,\lambda}(\Lambda^4T^*\Nhat)).
\end{align*}
Recall that $\hat{\mathcal{T}}\subseteq\R^n\cong
T_0\hat{\mathcal{T}}$ is open.  Then $X_2$, $Y_2$ and $Z_2$ are
Banach spaces, $U_2$ is an open neighbourhood of $(0,0,0)$ in $X_2$
and hence $W_2=U_2\times Y_2$ is an open neighbourhood of
$(0,0,0,0)$ in $X_2\times Y_2$.  Define
$\mathcal{G}_2:W_2\rightarrow Z_2$ by:
$$\mathcal{G}_2(\alpha,\beta,t,\gamma)=G_2(\alpha,\beta,t)+\gamma.$$
Then $d\mathcal{G}_2|_{(0,0,0,0)}:X_2\times Y_2\rightarrow Z_2$ acts
as
$$(\alpha,\beta,t,\gamma)\longmapsto d\alpha+d^*\beta+L_2(t)+\gamma.$$
By construction, $d\mathcal{G}_2|_{(0,0,0,0)}$ is surjective and,
using the fact that the image of $dG_2|_{(0,0,0)}$ meets $Y_2$ at
$0$ only,
\begin{align*}
\text{Ker}\,d\mathcal{G}_2|_{(0,0,0,0)}&=\{(\alpha,\beta,t,\gamma)\in
X_2\times
Y_2\,:\,d\alpha+d^*\beta+L_2(t)+\gamma=0\}\\
&\cong\{(\alpha,\beta,t)\in
X_2\,:\,d\alpha+d^*\beta+L_2(t)=0\}\cong\mathcal{I}_2(N,\lambda).
\end{align*}

\noindent As for Problem 1, Theorem \ref{ch6s2thm2} gives us that
$\text{Ker}\,\mathcal{G}_2$ is a smooth manifold near zero which may
be identified with an open neighbourhood
$\hat{\mathcal{M}}_2(N,\lambda)$ of $(0,0)$ in
$\mathcal{I}_2(N,\lambda)$.  We can again define a smooth map
$\pi_2:\hat{\mathcal{M}}_2(N,\lambda)\rightarrow\mathcal{O}_2(N,\lambda)$
such that $\text{Ker}\,\pi_2$ is locally homeomorphic near $(0,0)$
to an open neighbourhood of $N$ in $\mathcal{M}_2(N,\lambda)$.  We
thus have the following theorem.

\begin{thm}\label{ch8s3subs2thm1} Use the notation of Definitions
\ref{ch8s2subs2dfn1}, \ref{ch8s3subs2dfn1} and \ref{ch8s3subs2dfn2}.
There exists a smooth manifold $\hat{\mathcal{M}}_2(N,\lambda)$,
which is an open neighbourhood of $(0,0)$ in
$\mathcal{I}_2(N,\lambda)$, and a smooth map
$\pi_2:\hat{\mathcal{M}}_2(N,\lambda)\rightarrow\mathcal{O}_2(N,\lambda)$,
with $\pi_2(0,0)=0$, such that an open neighbourhood of zero in
$\text{\emph{Ker}}\,\pi_2$ is homeomorphic to an open neighbourhood
of $N$ in $\mathcal{M}_2(N,\lambda)$.
\end{thm}

\noindent We deduce that, if $\mathcal{O}_2(N,\lambda)=\{0\}$, the
moduli space for Problem 2 is a smooth manifold near $N$ of
dimension
$\text{dim}\,\mathcal{I}_2(N,\lambda)=\text{dim}\,\mathcal{I}_1(N,\lambda)
+\text{dim}\,\hat{\mathcal{T}}$, which we expect to occur for
generic choices of $N$ and the torsion-free $\text{G}_2$ structure
on $M$. We shall see, in $\S$\ref{ch8s5}, that if we choose a
suitable generic \emph{closed} $\G2$ structure on $M$ we may drop
the assumption that $N$ is generic and still obtain a smooth moduli
space.

\subsection{Problem 3}\label{ch8s3subs3}

We presume in this subsection that the reader is sufficiently
familiar with the schemata we have used in the previous two
subsections to be able to generalise them to Problem 3.  This allows
us to present a tidier treatment of the problem.

Recall the notation of Definitions \ref{ch8s2subs3dfn1},
\ref{ch8s2subs3dfn3} and \ref{ch8s2subs3dfn4}.

\begin{dfn}\label{ch8s3subs3dfn1} The \emph{infinitesimal deformation space}
$\mathcal{I}_3(N,\lambda)$ for Problem 3 is
\begin{align*}
\mathcal{I}_3(N,\lambda)&=\{(\alpha,t,f)\in
C^{\infty}_\lambda(\Lambda^2_+T^*\Nhat)\oplus
T_0\hat{\mathcal{T}}\oplus T_0\hat{\mathcal{F}}\,:\,d\alpha+L_2(t)+L_3(f)=0\}\\
&\cong\{(\alpha,\beta,t,f)\in
C^{\infty}_\lambda(\Lambda^2_+T^*\Nhat\oplus
\Lambda^4T^*\Nhat)\oplus T_0\hat{\mathcal{T}}\oplus
T_0\hat{\mathcal{F}}\\&\qquad\qquad\qquad\!\!\!\hspace{-0.5pt}\qquad\qquad\qquad\qquad\qquad\,:\,
d\alpha+d^*\beta+L_2(t)+L_3(f)=0\}.
\end{align*}
By Proposition \ref{ch8s2subs3prop5},
\begin{align*}\mathcal{I}_3(N,\lambda)&\cong\{(\alpha,\beta,t,f)\in
L^p_{k+1,\,\lambda}(\Lambda^2_+T^*\Nhat\oplus\Lambda^4T^*\Nhat)\oplus
T_0\hat{\mathcal{T}}\oplus
T_0\hat{\mathcal{F}}\\&\qquad\qquad\qquad\qquad\qquad\qquad\qquad\,:\,
d\alpha+d^*\beta+L_2(t)+L_3(f)=0\}.\end{align*}
\end{dfn}

\vspace{-20pt}

In considering obstructions, we first have the generalisation of
Propositions \ref{ch8s3subs1prop1} and \ref{ch8s3subs2prop1}.

\begin{prop}\label{ch8s3subs3prop1}
$F_3\big(L^p_{k+1,\,\lambda}(\Uhat)\times\hat{\mathcal{T}}\times\hat{\mathcal{F}}\big)\subseteq
d(L^p_{k+1,\,\lambda}(\Lambda^2T^*\Nhat))$.
\end{prop}

\noindent The proposition is proved in a similar way to Proposition
\ref{ch8s3subs2prop1} and so we omit the details.  The result leads
us to define our final obstruction space.

\begin{dfn}\label{ch8s3subs3dfn2}
From Propositions \ref{ch8s3subs1prop2} and \ref{ch8s3subs3prop1},
since $L_2$ and $L_3$ are linear maps on finite-dimensional vector
spaces, there exists a finite-dimensional subspace
$\mathcal{O}_3(N,\lambda)$ of
$L^p_{k,\,\lambda-1}(\Lambda^3T^*\Nhat)$ such that
$$
d(L^p_{k+1,\,\lambda}(\Lambda^2T^*\Nhat))=(d(L^p_{k+1,\,\lambda}(\Lambda^2_+T^*\Nhat))+
L_2(T_0\hat{\mathcal{T}})+L_3(
T_0\hat{\mathcal{F}}))\oplus\mathcal{O}_3(N,\lambda).$$ We define
$\mathcal{O}_3(N,\lambda)$ to be the \emph{obstruction space} for
Problem 3.
\end{dfn}

\begin{note} We may choose our obstruction spaces such that
$\mathcal{O}_3(N,\lambda)\subseteq\mathcal{O}_2(N,\lambda)\subseteq\mathcal{O}_1(N,\lambda)$.\end{note}

The use of the Implicit Function Theorem (Theorem \ref{ch6s2thm2})
in the derivation of Theorems \ref{ch8s3subs1thm1} and
\ref{ch8s3subs2thm1} can be easily generalised to give the
following.

\begin{thm}\label{ch8s3subs3thm1} Use the notation of Definitions \ref{ch8s2subs3dfn2},
\ref{ch8s3subs3dfn1} and \ref{ch8s3subs3dfn2}. There exists a smooth
manifold $\hat{\mathcal{M}}_3(N,\lambda)$, which is an open
neighbourhood of $(0,0,0)$ in $\mathcal{I}_3(N,\lambda)$, and a
smooth map
$\pi_3:\hat{\mathcal{M}}_3(N,\lambda)\rightarrow\mathcal{O}_3(N,\lambda)$,
with $\pi_3(0,0,0)=0$, such that an open neighbourhood of zero in
$\text{\emph{Ker}}\,\pi_3$ is homeomorphic to an open neighbourhood
of $(N,0)$ in $\mathcal{M}_3(N,\lambda)$.  
\end{thm}

\noindent We deduce that, if $\mathcal{O}_3(N,\lambda)=\{0\}$,
$\mathcal{M}_3(N,\lambda)$ is a smooth manifold near $(N,0)$ of
dimension
$\text{dim}\,\mathcal{I}_3(N,\lambda)=\text{dim}\,\mathcal{I}_2(N,\lambda)
+\text{dim}\,\hat{\mathcal{F}}$.  Moreover, the projection map
$\pi_{\hat{\mathcal{F}}}:\mathcal{M}_3(N,\lambda)\rightarrow\hat{\mathcal{F}}$
is smooth near $(N,0)$. We expect this to occur for generic choices
of $N$ and the torsion-free $\text{G}_2$ structure on $M$.  If we
allow ourselves to work with \emph{closed} $\G2$ structures on $M$,
we shall show in $\S$\ref{ch8s5} that we may drop our genericity
assumptions for $N$ and $(\varphi,g)$ and still get a smooth moduli
space.

\section{Dimension calculations}\label{ch8s4}

We shall relate the expected dimension of the moduli space for
Problem 1 to the index of $d+d^*$ as discussed in
$\S$\ref{indextheory}. Recall that $p>4$, $k\geq 2$ and
$\lambda\in(1,2)\setminus\mathcal{D}$.

\begin{dfn} Define
\begin{equation*}
\mathcal{H}^m=\{\xi\in L^2(\Lambda^mT^*\Nhat)\,:\,d\xi=d^*\xi=0\}.
\end{equation*}
 The Hodge star maps $\mathcal{H}^2$ into itself, so there is a
splitting $\mathcal{H}^2=\mathcal{H}_+^2\oplus\mathcal{H}_-^2$ where
$$\mathcal{H}_{\pm}^2=\mathcal{H}^2\cap
C^\infty(\Lambda_{\pm}^2T^*\Nhat).$$

Let $\mathcal{J}=\jmath\left(H_{\text{cs}}^2(\Nhat)\right)$, where
$H_\text{cs}^m(\Nhat)$ is the $m$th \emph{compactly supported}
cohomology group on $\Nhat$ and
$\jmath:H_\text{cs}^2(\Nhat)\rightarrow H^2_\text{dR}(\Nhat)$ is the
inclusion map. If $\alpha,\beta\in \mathcal{J}$, there exist
compactly supported closed 2-forms $\xi$ and $\eta$ such that
$\alpha=[\xi]$ and $\beta=[\eta]$. We define a product on
$\mathcal{J}\times \mathcal{J}$ by
\begin{equation}\label{ch6s5eq2}
\alpha\cup\beta=\int_{\Nhat} \xi\w\eta. \end{equation} Suppose that
$\xi^{\prime}$ and $\eta^{\prime}$ are also compactly supported with
$\alpha=[\xi^{\prime}]$ and $\beta=[\eta^{\prime}]$.  Then there
exist 1-forms $\chi$ and $\zeta$ such that $\xi-\xi^{\prime}=d\chi$
and $\eta-\eta^{\prime}=d\zeta$. Therefore,
\begin{align*}
\int_{\Nhat}\xi^{\prime}\w\eta^{\prime}&=\int_{\Nhat}(\xi-d\chi)\w(\eta-d\zeta)
=\int_{\Nhat}\xi\w\eta-d\chi\w\eta-\xi^{\prime}\w d\zeta\\[4pt]
&=\int_{\Nhat}\xi\w\eta-d(\chi\w\eta)-d(\xi^{\prime}\w\zeta)=\int_{\Nhat}\xi\w\eta,
\end{align*}
as both $\chi\w\eta$ and $\xi^{\prime}\w\zeta$ have compact support.
The product \eq{ch6s5eq2} on $\mathcal{J}\times \mathcal{J}$ is thus
well-defined and is a symmetric topological product with a signature
$(a,b)$.  By \cite[Example (0.16)]{Lockhart}, $\mathcal{H}^2\cong
\mathcal{J}$ and the isomorphism is given by $\xi\mapsto[\xi]$.
Thus, $\text{dim}\,\mathcal{H}_+^2=a$ and hence is a topological
number.
\end{dfn}

\begin{dfn}\label{ch8s4dfn1} Let
\begin{equation*}
(d_++d^*)_\lambda=d+d^*:L^p_{k+1,\,\lambda}(\Lambda^2_+T^*\Nhat\oplus\Lambda^4T^*\Nhat)\rightarrow
L^p_{k,\,\lambda-1}(\Lambda^3T^*\Nhat).
\end{equation*}
By Definition \ref{ch8s3subs1dfn1}, $\mathcal{I}_1(N,\lambda)$ is
isomorphic to the kernel of this map.  Define the adjoint map by
\begin{equation*}
(d^*_++d)_{-3-\lambda}=d^*_++d:L^q_{l+1,\,-3-\lambda}(\Lambda^3T^*\Nhat)\longrightarrow
L^q_{l,\,-4-\lambda}(\Lambda^2_+T^*\Nhat\oplus\Lambda^4T^*\Nhat),
\end{equation*}
where $q>1$ such that $1/p+1/q=1$ and $l\geq 4$.  The cokernel of
$(d_++d^*)_\lambda$ is then isomorphic to the kernel of
$(d^*_++d)_{-3-\lambda}$.
\end{dfn}

\begin{note}
The choice of $l\geq 4$ in Definition \ref{ch8s4dfn1} ensures that
$L^q_{l+1,\,-3-\lambda}\hookrightarrow C^{1,\,a}_{-3-\lambda}$ for
$0<a\leq 4-\frac{4}{q}=\frac{4}{p}<1$ by Theorem \ref{wembed}.
\end{note}

We now study the dimension of the kernel and cokernel of
$(d_++d^*)_{\mu}$.

\begin{prop}\label{ch8s4prop1} The kernel of $(d_++d^*)_{-2}$ is
isomorphic to $\mathcal{H}^2_+$.  Furthermore, if $\mu>-2$ is such
that $(-2,\mu]\cap\mathcal{D}=\emptyset$,
$\text{\emph{dim}}\,\text{\emph{Ker}}\,(d_++d^*)_{\mu}=\text{\emph{dim}}\,\mathcal{H}^2_+$.
\end{prop}
\begin{proof}
Using \eq{ch6s2eq3} and the Maximum Principle,
\begin{align*}\mathcal{H}^2_+&=\{\alpha\in L^2(\Lambda^2T^*\Nhat)\cap
C^{\infty}(\Lambda^2_+T^*\Nhat) \,:\,d\alpha=d^*\alpha=0\}\\
&=\{\alpha\in L^2_{0,\,-2}(\Lambda^2_+T^*\Nhat)\cap
C^{\infty}(\Lambda^2_+T^*\Nhat)\,:\,d\alpha=0\}\\
&\cong\{(\alpha,\beta)\in L^2_{0,\,-2}(\Lambda^2_+T^*\Nhat\oplus
\Lambda^4T^*\Nhat)\,:\,\alpha\in C^\infty(\Lambda^2_+T^*\Nhat),\;
d\alpha+d^*\beta=0\}.\end{align*} This gives the first part of the
proposition.


If $-2\notin\mathcal{D}$, $[-2,\mu]\cap\mathcal{D}=\emptyset$ and
thus, by Proposition \ref{nochange},
$\text{dim}\,\text{Ker}\,(d_++d^*)_{\mu}=\text{dim}\,\text{Ker}\,(d_++d^*)_{-2}$.

Suppose now that $-2\in\mathcal{D}$ and that $(\alpha,\beta)$
corresponds to a self-dual 2-form and 4-form on $\Nhat$ which are
subtracted from the kernel of $(d_++d^*)_{\nu}$ as $\nu$ crosses
$-2$ from below.  By the work in \cite[$\S$3 \&
$\S$4]{LockhartMcOwen} this occurs if and only if $(\alpha,\beta)$
is asymptotic to an $O(r^{-2})$ form $\xi$ on $C_i$, for some $i$,
satisfying $(d+d^*)\xi=0$.  (The form $\xi$ is determined by an
element of $D(-2,i)$, using the notation of Proposition
\ref{ch6s3prop1}.) Therefore, $(\alpha,\beta)$ is of order
$O(\rho^{-2})$ as $\rho\rightarrow 0$ and thus lies in $L^2$.  We
deduce that $(\alpha,\beta)\in\,\text{Ker}\,(d_++d^*)_{-2}$,
implying that the function $k(\nu)=\,\text{Ker}\,(d_++d^*)_{\nu}$ is
upper semi-continuous at $-2$ by Proposition \ref{nochange}.

The second part of the proposition is thus proved.
\end{proof}

\begin{prop}\label{ch8s4prop3} If $\mu<-1$ is such that
$[\mu,-1)\cap\mathcal{D}=\emptyset$, the cokernel of
$(d_++d^*)_{\mu}$ is isomorphic to $H^1_{\text{\emph{dR}}}(\Nhat)$.
\end{prop}

\begin{proof}
By Theorem \ref{ch6s3thm1}, there exists a countable discrete subset
$\mathcal{D}^\prime$ of rates $\nu$ such that
\begin{equation}\label{ch7s2eq3}
d+d^*:L_{k+1,\,\nu}^{p}(\Lambda^\text{even}T^*\Nhat) \rightarrow
L_{k,\,\nu-1}^{p}(\Lambda^\text{odd}T^*\Nhat) \end{equation} is not
Fredholm. Clearly, $\mathcal{D}^\prime\supseteq\mathcal{D}$. For
$\nu\notin\mathcal{D}^\prime$ with $\nu<-1$, so that $-3-\nu>\nu-1$,
$$
L_{k,\,\nu-1}^{p}(\Lambda^\text{odd}T^*\Nhat)=
(d+d^*)(L_{k+1,\,\nu}^{p}(\Lambda^\text{even}T^*\Nhat))
\oplus\mathcal{K},
$$ where $\mathcal{K}$ is the kernel of the adjoint map
\begin{equation*}
d+d^*:L_{l+1,\,-\nu-3}^{q}(\Lambda^\text{odd}T^*\Nhat)\rightarrow
L_{l,\,-\nu-4}^{q}(\Lambda^\text{even}T^*\Nhat),
\end{equation*}
for $1/p+1/q=1$ and $l\geq 4$, which is graded and closed under the
Hodge star.

If $\gamma\in L_{k,\,\nu-1}^{p}(\Lambda^3T^*\Nhat)$ then
$(\ast\gamma,\gamma)\in
L_{k,\,\nu-1}^{p}(\Lambda^\text{odd}T^*\Nhat)$ and hence there exist
some $\gamma_m\in L_{k+1,\,\nu}^{p}(\Lambda^mT^*\Nhat)$, for
$m=0,2,4$, and $\eta\in\mathcal{K}$ such that
\begin{equation*}
(\ast\gamma,\gamma)=
(d+d^*)(\gamma_0,\gamma_2,\gamma_4)+\eta.\end{equation*} By applying
the Hodge star,
\begin{equation*}
(\ast\gamma,\gamma)=
(d+d^*)(\ast\gamma_4,\ast\gamma_2,\ast\gamma_0)+\ast
\eta.\end{equation*} Adding the above formulae and averaging gives:
$$\gamma=d\!\left(\frac{\gamma_2+\ast\gamma_2}{2}\right)+
d^*\!\left(\frac{\ast\gamma_0+\gamma_4}{2}\right)+\tilde{\eta}$$
where $\tilde{\eta}\in\mathcal{K}\cap
L_{k,\,\nu-1}^{p}(\Lambda^3T^*\Nhat)$.  We deduce that
$$L_{k,\,\nu-1}^{p}(\Lambda^3T^*\Nhat)=\Big(d(L_{k+1,\,\nu}^{p}(\Lambda_+^2T^*\Nhat))
+d^*(L_{k+1,\,\nu}^{p}(\Lambda^4T^*\Nhat))\Big)\oplus\mathcal{K}^3,$$
where $\mathcal{K}^3=\mathcal{K}\cap
L_{k,\,\nu-1}^{p}(\Lambda^3T^*\Nhat)$.  Moreover, for
$\nu\notin\mathcal{D}^\prime$, $\nu<-1$,
\begin{equation*}
d(L_{k+1,\,\lambda}^{p}(\Lambda^2T^*\Nhat))=d(L_{k+1,\,\lambda}^{p}(\Lambda_+^2T^*\Nhat)).
\end{equation*}
We must surely have that the images are equal for
$\nu\notin\mathcal{D}$, $\nu<-1$, as well.

Thus, the cokernel of $(d_++d^*)_{\mu}$ is isomorphic to the kernel
of
\begin{equation}\label{ch8s4eq5}
(d^*+d)_{-3-\mu}=d^*+d:L^q_{l+1,\,-3-\mu}(\Lambda^3T^*\Nhat)\longrightarrow
L^q_{l,\,-4-\mu}(\Lambda^2T^*\Nhat\oplus\Lambda^4T^*\Nhat).
\end{equation}
Using \eq{ch6s2eq3} as in the proof of Proposition \ref{ch8s4prop1},
the kernel of $(d^*+d)_{-3-(-1)}=(d^*+d)_{-2}$ is isomorphic to
$\mathcal{H}^3$.
  By \cite[Example (0.16)]{Lockhart}, $\mathcal{H}^3\cong H^1_{\text{dR}}(\Nhat)$ and
  the isomorphism is given by $\gamma\mapsto[*\gamma]$.  Since
$[\mu,-1)\cap\mathcal{D}=\emptyset$, there are no changes in the
cokernel in $[\mu,-1)$ by Proposition \ref{nochange}. Moreover, the
dimension of the cokernel is lower semi-continuous in $\mu$ at $-1$;
this fact can be demonstrated using similar methods to those
employed in the proof of Proposition \ref{ch8s4prop1}. The result
follows.
\end{proof}




By Proposition \ref{nopoints}, $(-2,-1]\cap\mathcal{D}=\emptyset$.
Therefore, for any $\mu\in(-2,-1]$,
$$\text{dim}\,\text{Ker}\,(d_++d^*)_{\mu}=\text{dim}\,\mathcal{H}^2_+
\quad\text{and}\quad\text{dim}\,\text{Coker}\,(d_++d^*)_{\mu}=
b^1(\Nhat),$$ using Propositions \ref{ch8s4prop1} and
\ref{ch8s4prop3}.  Knowing the index of $(d_++d^*)_{\mu}$ for
$\mu\in(-2,-1]$, we can calculate it for all growth rates using
Theorem \ref{ch6s3thm2}.

\begin{prop}\label{ch8s4prop6}
Use the notation of Propositions \ref{ch6s3prop1} and
\ref{ch6s3prop3}. If $\lambda\in(1,2)$, $\lambda\notin\mathcal{D}$,
the index of $(d_++d^*)_\lambda$ is given by:
$$\text{\emph{ind}}\,(d_++d^*)_\lambda=\text{\emph{dim}}\,\mathcal{H}^2_+-b^1(\Nhat)
\,\,-\!\!\!\!\!\!\sum_{\mu\in(-1,\lambda)\cap\mathcal{D}}\!\!\!\!\!\d(\mu).$$
\end{prop}

However, the obstruction space $\mathcal{O}_1(N,\lambda)$ given in
Definition \ref{ch8s3subs1dfn2} is a subspace of the cokernel of
$(d_++d^*)_\lambda$, so we must relate their dimensions.

\begin{prop}\label{ch8s4prop5}
The following inequality holds:
$$\text{\emph{dim}}\,\mathcal{O}_1(N,\lambda)\leq\text{\emph{dim}}\,
\text{\emph{Coker}}\,(d_++d^*)_{\lambda}-b^1(\Nhat).$$
\end{prop}

\begin{proof}
From the proof of Proposition \ref{ch8s3subs1prop1}, the image of
\begin{equation*}
(d+d^*)_\lambda=d+d^*:L^p_{k+1,\,\lambda}(\Lambda^2T^*\Nhat\oplus\Lambda^4T^*\Nhat)
\rightarrow L^p_{k,\,\lambda-1}(\Lambda^3T^*\Nhat)\end{equation*} is
characterised as the subspace of
$L^p_{k,\,\lambda-1}(\Lambda^3T^*\Nhat)$ which is $L^2$-orthogonal
to the kernel $\mathcal{K}$ of $(d^*+d)_{-3-\lambda}$ defined by
\eq{ch8s4eq5}.
Furthermore, as noticed in the proof of Proposition
\ref{ch8s3subs1prop2}, $\text{Image}\,(d+d^*)_{\lambda}$ has finite
codimension in $L^p_{k,\,\lambda-1}(\Lambda^3T^*\Nhat)$.  Therefore,
we may choose a finite-dimensional space $\mathcal{C}$ of smooth
compactly supported 3-forms on $\Nhat$ such that
$$L^p_{k,\,\lambda-1}(\Lambda^3T^*\Nhat)=\,\text{Image}\,(d+d^*)_\lambda\oplus\mathcal{C}$$
and so that the product $\mathcal{C}\times\mathcal{K}\rightarrow\R$
given by $(\gamma,\eta)\mapsto\langle\gamma,\eta\rangle_{L^2}$ is
nondegenerate.

We may similarly deduce that the image of $(d_++d^*)_\lambda$ is the
subspace of $L^p_{k,\,\lambda-1}(\Lambda^3T^*\Nhat)$ which is
$L^2$-orthogonal to the kernel $\mathcal{K}^\prime$ of
$(d^*_++d)_{-3-\lambda}$. Then
$\mathcal{K}^\prime\supseteq\mathcal{K}$ and $\mathcal{K}$ consists
of closed and coclosed 3-forms, whereas $\mathcal{K}^\prime$
consists of 3-forms $\eta$ such that $d\eta=d^*_+\eta=0$. Hence, we
may choose a subspace $\mathcal{K}^{\prime\prime}$ of
$\mathcal{K}^\prime$, transverse to $\mathcal{K}$, comprising
3-forms which are \emph{not} coclosed and such that
$\mathcal{K}^\prime=\mathcal{K}\oplus\mathcal{K}^{\prime\prime}$.

The next stage is to extend $\mathcal{C}$ to a space
$\mathcal{C}^\prime=\mathcal{C}\oplus\mathcal{C}^{\prime\prime}$,
where $\mathcal{C}^{\prime\prime}$ consists of smooth \emph{exact}
compactly supported 3-forms on $\Nhat$, such that
$$L^p_{k,\,\lambda-1}(\Lambda^3T^*\Nhat)=\,\text{Image}\,(d_++d^*)_\lambda\oplus\mathcal{C}^{\prime}$$
and such that the product
$\mathcal{C}^{\prime\prime}\times\mathcal{K}^{\prime\prime}\rightarrow\R$
given by $(\gamma,\eta)\mapsto\langle\gamma,\eta\rangle_{L^2}$ is
nondegenerate, which is possible as $\mathcal{K}^{\prime\prime}$
comprises forms which are not coclosed.  By construction,
$\mathcal{C}^{\prime\prime}$ is a valid choice for
$\mathcal{O}_1(N,\lambda)$ by Proposition \ref{ch8s3subs1prop2}.
Therefore,
$$\text{dim}\,\mathcal{O}_1(N,\lambda)=\,
\text{dim}\,\mathcal{C}^\prime-\text{dim}\,\mathcal{C}=\, \text{dim
Coker}\,(d_++d^*)_{\lambda}-\text{dim}\,\mathcal{K}.$$ If $\gamma$
lies in the kernel of \eq{ch8s4eq5} for rate $\mu=-1$ then
$\gamma\in\mathcal{K}$ for $\lambda\in(1,2)$
 by Theorem
\ref{ch6s2thm1}. Thus, the map from $\mathcal{K}$ to
$H^1_{\text{dR}}(\Nhat)$ given by $\gamma\mapsto[*\gamma]$ is
surjective. This gives the result.
\end{proof}

We may now calculate a lower bound for the expected dimension of
$\mathcal{M}_1(N,\lambda)$ using Propositions \ref{ch8s4prop6} and
\ref{ch8s4prop5}.

\begin{prop}\label{ch8s4prop7}
Using the notation of Propositions \ref{ch6s3prop1} and
\ref{ch6s3prop3},
$$\text{\emph{dim}}\,\mathcal{I}_1(N,\lambda)-\text{\emph{dim}}\,\mathcal{O}_1(N,\lambda)
\geq\text{\emph{dim}}\,\mathcal{H}^2_+
\,-\!\!\!\!\!\!\sum_{\mu\in(-1,\lambda)\cap\mathcal{D}}\!\!\!\!\!\d(\mu).$$
\end{prop}

\noindent Recalling that the dimension of $\mathcal{T}$ given in
Definition \ref{ch8s2subs2dfn2} is $21s$, we derive analogous
results for our other problems.

\begin{prop}\label{ch8s4prop8} Using
the notation of Definitions \ref{ch8s2subs2dfn2} and
\ref{ch8s2subs3dfn1} and Propositions \ref{ch6s3prop1} and
\ref{ch6s3prop3},
\begin{align*}
\text{\emph{dim}}\,\mathcal{I}_2(N,\lambda)-\text{\emph{dim}}\,\mathcal{O}_2(N,\lambda)
&\geq\text{\emph{dim}}\,\mathcal{H}^2_++21s-\text{\emph{dim$\,$H}}
-\!\!\!\!\!\sum_{\mu\in(-1,\lambda)\cap\mathcal{D}}\!\!\!\!\!\d(\mu).
\intertext{and}
\text{\emph{dim}}\,\mathcal{I}_3(N,\lambda)-\text{\emph{dim}}\,\mathcal{O}_3(N,\lambda)
&\geq\text{\emph{dim}}\,\mathcal{H}^2_++21s-\text{\emph{dim$\,$H}}+\text{\emph{dim}}\,\hat{\mathcal{F}}
-\!\!\!\!\!\sum_{\mu\in(-1,\lambda)\cap\mathcal{D}}\!\!\!\!\!\d(\mu).
\end{align*}
\end{prop}

We note that Propositions \ref{nopoints}, \ref{ch8s4prop3},
\ref{ch8s4prop5} and \ref{ch8s4prop7} imply the following bound on
$\text{dim}\,\mathcal{O}_1(N,\lambda)$.

\begin{prop}\label{ch8s4prop9} In the notation of Propositions
\ref{ch6s3prop1} and \ref{ch6s3prop3},
$$\text{\emph{dim}}\,\mathcal{O}_1(N,\lambda)\leq\sum_{\mu\in(-1,\lambda)\cap\mathcal{D}}
\!\!\!\!\d(\mu)$$
\end{prop}

\noindent We also know that, in Problem 2, we remove the
obstructions which correspond to translations of the singularities
and $\G2$ transformations of the tangent cones.  These obstructions
occur, respectively, at rates $0$ and $1$.  Hence, $\d(0)\geq 7s$,
$\d(1)\geq 14s-\text{dim H}$ and we have the following stronger
bound on the dimension of $\mathcal{O}_2(N,\lambda)$.

\begin{prop} In the notation of Definition \ref{ch8s2subs2dfn2} and Propositions
\ref{ch6s3prop1} and \ref{ch6s3prop3},
$$\text{\emph{dim}}\,\mathcal{O}_2(N,\lambda)\leq-21s+\text{\emph{dim}}\,\text{\emph{H}}+
\!\!\!\!\!\sum_{\mu\in(-1,\lambda)\cap\mathcal{D}}
\!\!\!\!\d(\mu).$$
\end{prop}

\section{$\varphi$-Closed 7-manifolds}\label{ch8s5}

For our deformation problems we have assumed the ambient manifold
$(M,\varphi,g)$ is a $\text{G}_2$ manifold; that is, $M$ is endowed
with a $\G2$ structure such that $d\varphi=d^*\varphi=0$. However,
the results of McLean \cite{McLean} we have used, which are based
upon the linearisation of the map we denoted $F_1$ in Definition
\ref{ch8s2subs1dfn2}, still hold if this condition on $\varphi$ is
relaxed to just $d\varphi=0$. Thus, our deformation theory results
hold if $(M,\varphi,g)$ is a $\varphi$-closed 7-manifold in the
sense of Definition \ref{ch2s3subs2dfn3}. \begin{remark} The effect
of $*\varphi$ not being closed on $M$ means that coassociative
4-folds in $M$ are no longer necessarily volume minimizing in their
homology class. This does not, however, affect our
discussion.\end{remark}

The use of $\varphi$-closed 7-manifolds $(M,\varphi,g)$ is that
closed $\text{G}_2$ structures occur in infinite-dimensional
families, since the set of closed definite 3-forms on $M$, in the
sense of Definition \ref{ch2s3subs2dfn1}, is open. We show that we
can choose a family $\mathcal{F}$, in a similar fashion to
Definition \ref{ch8s2subs3dfn1} of Problem 3, of closed $\text{G}_2$
structures on $M$ such that
$\text{dim}\,\mathcal{F}=\text{dim}\,\mathcal{O}_1(N,\lambda)$ and,
further, such that $\mathcal{O}_3(N,\lambda)=\{0\}$. In other words,
we have enough freedom in our choice of $\mathcal{F}$ to ensure that
$dF_3|_{(0,0,0)}$, as given in Definition \ref{ch8s2subs3dfn3}, maps
onto $d(L^p_{k+1,\,\lambda}(\Lambda^2T^*\Nhat))$.  Then
$\mathcal{M}_3(N,\lambda)$ is a smooth manifold near $(N,0)$ by
Theorem \ref{ch8s3subs3thm1} and
$\pi_{\hat{\mathcal{F}}}:\mathcal{M}_3(N,\lambda)\rightarrow\hat{\mathcal{F}}$
is a smooth map near $(N,0)$.

\emph{Sard's Theorem} \cite[p. 173]{Lang} states that, if
$f:X\rightarrow Y$ is a smooth map between finite-dimensional
manifolds, the set of $y\in Y$ with some $x\in f^{-1}(y)$ such that
$df|_x:T_xX\rightarrow T_yY$ is \emph{not} surjective is of measure
zero in $Y$. Therefore, $f^{-1}(y)$ is a submanifold of $X$ for
almost all $y\in Y$.

By Sard's Theorem, $\pi_{\hat{\mathcal{F}}}^{-1}(f)$ is a smooth
manifold near $(N,f)$ for almost all $f\in\hat{\mathcal{F}}$. As
observed in Definition \ref{ch8s2subs3dfn2},
$\pi_{\hat{\mathcal{F}}}^{-1}(f)$ corresponds to the moduli space of
deformations for Problem 2 defined using the $\G2$ structure
$(\varphi^f,g^f)$. Thus, for any given $N$, a generic perturbation
of the closed $\text{G}_2$ structure within $\mathcal{F}$ ensures
that $\mathcal{M}_2(N,\lambda)$ is smooth near $N$.

We thus prove the following, which is similar to the result
\cite[Theorem 9.1]{Joyce2}.

\begin{thm} Let $(M,\varphi,g)$ be a $\varphi$-closed 7-manifold in
the sense of Definition \ref{ch2s3subs2dfn3} and let $N$ in
$(M,\varphi,g)$ be a CS coassociative 4-fold at $z_1,\ldots,z_s$
with rate $\lambda\in (1,2)\setminus\mathcal{D}$, where
$\mathcal{D}$ is defined in Proposition \ref{ch6s3prop1}.  Use the
notation of Definitions \ref{ch8s2subs2dfn2}, \ref{ch8s3subs1dfn2}
and \ref{ch8s3subs3dfn2} and Proposition \ref{ch6s3prop3}. Let
$m=\text{\emph{dim}}\,\mathcal{O}_1(N,\lambda)$ and let
$\hat{\mathcal{F}}$ be an open ball about $0$ in $\R^m$. There
exists a smooth family
$\mathcal{F}=\{(\varphi^f,g^f)\,:\,f\in\hat{\mathcal{F}}\}$ of
closed $\text{\emph{G}}_2$ structures on $M$ such that
$\mathcal{O}_3(N,\lambda)=\{0\}$. Hence, the moduli space of
deformations for Problem 3 is a smooth manifold near $(N,0)$ of
dimension greater than or equal to
$$\text{\emph{dim}}\,\mathcal{H}^2_++21s-\text{\emph{dim}}\,\text{\emph{H}}
+\text{\emph{dim}}\,\mathcal{O}_1(N,\lambda)
-\!\!\!\!\!\sum_{\mu\in(-1,\lambda)\cap\mathcal{D}}\!\!\!\!\!\d(\mu).$$
Moreover, for generic $f\in\hat{\mathcal{F}}$, the moduli space of
deformations in $(M,\varphi^f,g^f)$ for Problem 2 is a smooth
manifold near $N$ of dimension greater than or equal to
$$\text{\emph{dim}}\,\mathcal{H}^2_++21s-\text{\emph{dim}}\,\text{\emph{H}}
-\!\!\!\!\!\sum_{\mu\in(-1,\lambda)\cap\mathcal{D}}\!\!\!\!\!\d(\mu).$$
\end{thm}

\begin{proof}
Use the notation in the proof of Proposition \ref{ch8s4prop5}.
Recall that we have a subspace $\mathcal{K}^{\prime\prime}$ of
$L^q_{l+1,\,-3-\lambda}(\Lambda^3T^*\Nhat)$ consisting of forms
$\eta$ such that $d\eta=d^*_+\eta=0$ but $d^*\eta\neq0$.  Moreover,
$\mathcal{O}_1(N,\lambda)$ can be chosen to be a space of smooth
compactly supported exact 3-forms $\gamma$ such that
$\langle\gamma,\eta\rangle_{L^2}=0$ for all
$\eta\in\mathcal{K}^{\prime\prime}\setminus\{0\}$ implies that
$\gamma=0$.  Therefore
$\mathcal{K}^{\prime\prime}\cong(\mathcal{O}_1(N,\lambda))^*$ and
hence has dimension $m$.

Let $\{\eta_1,\ldots,\eta_m\}$ be a basis for
$\mathcal{K}^{\prime\prime}$ and choose a basis
$\{d\upsilon_1,\ldots,d\upsilon_m\}$ for $\mathcal{O}_1(N,\lambda)$,
where $\upsilon_j$ is a smooth compactly supported 2-form for all
$j$, such that $\langle
d\upsilon_i,\eta_j\rangle_{L^2}=\delta_{ij}$. This is possible
because the $L^2$ product on
$\mathcal{O}_1(N,\lambda)\times\mathcal{K}^{\prime\prime}$ is
nondegenerate.  For $f=(f_1,\ldots,f_m)\in\R^m$ define
$$\upsilon_f=\sum_{j=1}^mf_j\upsilon_j.$$
Using the notation of Proposition \ref{ch8s2subs3prop1} define
$(\varphi^f,g^f)$, for $f$ in a sufficiently small open ball
$\hat{\mathcal{F}}$ about $0$ in $\R^m$, to be a closed $\G2$
structure on $M$ such that $\Xi([\varphi^f|_T])=0$ in
$H^3_{\text{cs}}(\Nhat)$ and $\varphi^f|_{\Nhat}=d\upsilon_f$.
Recall from Definitions \ref{ch8s2subs3dfn3} and
\ref{ch8s3subs3dfn2}  that we have a linear map
$L_3:T_0\hat{\mathcal{F}}\cong\R^m\rightarrow
d(L^p_{k+1,\,\lambda}(\Lambda^2T^*\Nhat))$ arising from
$dF_3|_{(0,0,0)}$.  By construction, $L_3(f)=d\upsilon_f$ for
$f\in\R^m$ and hence $L_3$ maps onto $\mathcal{O}_1(N,\lambda)$.
Proposition \ref{ch8s3subs1prop2} and Definition
\ref{ch8s3subs3dfn2} imply that $\mathcal{O}_3(N,\lambda)=\{0\}$ as
required.

The latter parts of the theorem follow from the discussion preceding
it and Proposition \ref{ch8s4prop8}.
\end{proof}

\begin{ack}
Many thanks go to Dominic Joyce for his enormous help with this
project.  I would also like to thank Alexei Kovalev for useful
comments and suggestions and EPSRC for providing the funding for
this research.
\end{ack}



\begin{thebibliography}{99}

\bibitem{Bartnik} R. Bartnik, {\it The Mass of an Asymptotically
Flat Manifold}, Comm. Pure Appl. Math. {\bf 39}, 661-693, 1986.

\bibitem{Bryant4} R. L. Bryant, {\it Some Remarks on
$\text{\emph{G}}_2$-Structures}, preprint,\\
http://www.arxiv.org/math.DG/0305124, 2003


\bibitem{HarLaw} R. Harvey and H. B. Lawson, {\it Calibrated Geometries},
Acta Math. {\bf 148}, 47-152, 1982.

\bibitem{Federer} H. Federer, {\it Geometric Measure Theory},
Grundlehren Math. Wiss. \textbf{153}, Springer--Verlag, Berlin,
1969.

\bibitem{Joy1} D. D. Joyce, {\it Compact Manifolds with Special Holonomy},
OUP, Oxford, 2000.

\bibitem{Joyce1} D. D. Joyce, {\it Special Lagrangian Submanifolds
with Isolated Conical Singularities. I. Regularity}, Ann. Global
Ann. Geom. {\bf 25}, 201-251, 2004.
\bibitem{Joyce2} D. D. Joyce, {\it Special Lagrangian Submanifolds
with Isolated Conical Singularities. II. Moduli Spaces}, Ann. Global
Ann. Geom. {\bf 25}, 301-352, 2004.
\bibitem{Joyce3} D. D. Joyce, {\it Special Lagrangian Submanifolds
with Isolated Conical Singularities. III. Desingularization, The
Unobstructed Case}, Ann. Global Ann. Geom. {\bf 26}, 1-58, 2004.
\bibitem{Joyce4} D. D. Joyce, {\it Special Lagrangian Submanifolds
with Isolated Conical Singularities. IV. Desingularization,
Obstructions and Families}, Ann. Global Ann. Geom. {\bf 26},
117-174, 2004.
\bibitem{Joyce5} D. D. Joyce, {\it Special Lagrangian Submanifolds
with Isolated Conical Singularities. V. Survey and Applications}, J.
Differential Geom. {\bf 63}, 299-347, 2003.

\bibitem{Lang} S. Lang, {\it Differentiable Manifolds},
Addison--Wesley, Reading,\\ Massachusetts, 1972.
\bibitem{Lang2} S. Lang, {\it Real Analysis}, Second Edition, Addison--Wesley, Reading, Massachusetts, 1983.


\bibitem{Lockhart} R. B. Lockhart, {\it Fredholm, Hodge and Liouville
Theorems on Noncompact Manifolds}, Trans. Amer. Math. Soc. {\bf
301}, 1-35, 1987.
\bibitem{LockhartMcOwen} R. B. Lockhart and R. C. McOwen, {\it
Elliptic Differential Operators on Noncompact Manifolds}, Ann. Sc.
Norm. Super. Pisa Cl. Sci. {\bf 12}, 409-447, 1985.

\bibitem{Lotay} J. Lotay, \textit{Deformation Theory of
Asymptotically Conical Coassociative 4-folds}, preprint,
http://www.arxiv.org/math.DG/0411116, 2004.

\bibitem{McLean} R. C. McLean, {\it Deformations of Calibrated
Submanifolds}, Comm. Anal. Geom. {\bf 6}, 705-747, 1998.

\bibitem{Marshall} S. P. Marshall, {\it Deformations of Special
Lagrangian Submanifolds}, DPhil thesis, Oxford University, 2002.

\bibitem{Morrey} C. B. Morrey, {\it Multiple Integrals in the
Calculus of Variations}, Grundlehren Series Volume 130,
Springer--Verlag, Berlin, 1966.

\bibitem{Salamon} S. Salamon, {\it Riemannian Geometry and Holonomy
Groups}, Pitman Research Notes in Mathematics \textbf{201}, Longman,
Harlow, 1989.

\end{thebibliography}
\end{document}